\documentclass[10pt,reqno]{amsart}
\usepackage[margin=1.5 in]{geometry} %changes margins
\usepackage{amsfonts,amsmath, amssymb,mathrsfs,mathtools} %math fonts
\usepackage{enumitem,hyperref,cleveref} %for references, hyperlinks and enummerations
\usepackage{todonotes} %to insert comments

\usepackage[doi=false, url=false, eprint=false, maxbibnames=10, firstinits=true, isbn=false]{biblatex}
\DeclareUnicodeCharacter{0301}{*************************************}
\newtheorem{theorem}{Theorem}[section]
\newtheorem{proposition}[theorem]{Proposition}
\newtheorem{lemma}[theorem]{Lemma}
\newtheorem{corollary}[theorem]{Corollary}
\theoremstyle{definition}
\newtheorem{definition}[theorem]{Definition}
\theoremstyle{remark}
\newtheorem{remark}[theorem]{Remark}
\numberwithin{equation}{section}
\newtheorem{question}{Question}

\newcommand{\R}{{\mathbb R}}

\newcommand{\s}{{\mathbb S}}
\newcommand{\T}{{\mathbb T}}
\newcommand{\Z}{{\mathbb Z}}

\newcommand{\mcg}{_{\operatorname{McG}}}

\newcommand{\norm}[1]{\left\lVert#1\right\rVert}

\newcommand{\bb}[1]{\prescript{b}{}{#1}}
\newcommand{\btm}{\bb{T^\ast\!M}}
\newcommand{\bm}[1]{\prescript{b^m}{}{#1}}
\newcommand{\bmtm}{\bm{T^\ast\!M}}

\newcommand{\m}[1]{{\color{magenta}#1\color{black}}}

\makeatletter
\def\namedlabel#1#2{\begingroup
   \def\@currentlabel{#2}%
   \label{#1}\endgroup
}
\makeatother

\title[$b$-Contact structure on tentacular hyperboloids]{$b$-Contact structures on tentacular hyperboloids}

\author{Michael Vogel}
\address{Michael Vogel, Fakult\"{a}t f\"{u}r Mathematik, Ruhr-Universit\"{a}t Bochum, Bochum, Germany}
\curraddr{}
\email{michael.vogel-j6m@rub.de}

\author{Jagna Wi\'{s}niewska}

\date{\today}

\begin{document}

\maketitle

\begin{abstract}
This paper connects two different approaches to the analysis of Hamiltonian dynamics on non-compact energy hypersurfaces - $b$-symplectic geometry with its singular symplectic form and Floer techniques for tentacular Hamiltonians. More precisely, we show how to equip tentacular hyperboloids with a $b$-contact structure.
We construct a $b^3$-symplectic manifold $(X, Z, \omega_b)$, such that each connected component of $X\setminus Z$ is symplectomorphic to the standard symplectic space $(T^*\R^n, \omega_0)$. For a tentacular hyperboloid $S \subseteq T^*\R^n$ we look at its copies in $X \setminus Z$ and show that their completion in $(X,Z, \omega_b)$ is a smooth hypersurface of $b$-contact type.
\end{abstract}

\section{Introduction}

The question of existence of periodic orbits on a given energy level set is natural from the perspective of physical applications and has been a driving question in Hamiltonian dynamics and symplectic topology. The first general existence results for compact hypersurfaces were due to Rabinowitz \cite{Rabinowitz1979, Rabinowitz1978} and Weinstein \cite{weinstein1979} for the cases where the hypersurface is star-shaped and convex, respectively. Weinstein introduced the more general definition of a contact-type hypersurface, and used it to formulate the famous \emph{Weinstein Conjecture}, which asserts that if $S$ is a compact contact-type hypersurface in an arbitrary symplectic manifold $M$, then $S$ carries at least one periodic orbit. Later %, in 1987, 
Viterbo \cite{Viterbo1987} proved the conjecture for $M = \R^{2n}$.
On the other hand, ever since the introduction of Gromov’s holomorphic curves and Floer’s techniques \cite{Floer1989} to the symplectic world, the approach to this problem has often been based on the properties of spaces of holomorphic curves (or Floer trajectories) with prescribed asymptotics. For instance, Hofer’s proof of the Weinstein Conjecture for the three-dimensional sphere in \cite{Hofer1993} and Taubes’ proof for arbitrary closed, contact three-manifolds in \cite{Taubes2007} are obtained by studying holomorphic curves in the symplectization of the given contact manifold. However, the conjecture is still open for compact, contact-type hypersurfaces in arbitrary manifolds of dimension $2n\geq 4$.

The question of existence of periodic orbits in a non-compact setting is even more interesting from the perspective of applications, since it appears in many natural dynamical systems, like for example in celestial mechanics in particular the three body problem. However, up to recently little has been known about periodic orbits on non-compact energy hypersurfaces: even the Floer-theoretical methods break down if the compactness assumption is dropped. Additional geometric and topological assumptions are needed to make up for the lack of compactness. Nevertheless, in the recent years there have been new developments in the area of non-compact Hamiltonian and Reeb dynamics. In 2009 the Weinstein conjecture was answered for mechanical hypersurfaces in standard symplectic space and cotangent bundles in \cite{Pasquotto2009} and its sequel \cite{rot2016}, and in \cite{Suhr2016}.

In 2017 the first results in applying the Floer theoretical methods to non-compact energy hypersurfaces were presented by the second author in her PhD thesis \cite{Wisniewska2017} where she introduced a class of so called \emph{tentacular Hamiltonians} on $\R^{2n}$ with non-compact level sets for which Rabinowitz Floer homology is well defined. The non-compact level sets of tentacular Hamiltonians include a wide class of hyperboloids, called \emph{tentacular hyperboloids}. Later on, in 2020, Fauck, Merry and the second author \cite{Fauck2021} computed the Rabinowitz Floer homology for tentacular hyperboloids and used it to prove the Weinstein Conjeture for perturbations thereof.

Meanwhile, in 2017, Cieliebak, Eliashberg and Polterovich \cite{cieliebak2017} defined and computed symplectic homology for a subclass of tentacular hyperboloids, whereas Del Pino and Presas \cite{DelPino2017} proved the Weinstein conjecture for a special class of open contact manifolds, namely those that occur as overtwisted leaves of a contact foliation in a closed manifold. In 2020 Ganatra, Pardon and Shende \cite{ganatra2020} tackled the issue of non-compactness %of a contact-type hypersurface 
by introducing a class of Liouville manifolds with boundary called Liouville sectors for which they defined wrapped Fukaya cathegory and symplectic homology.

Another approach how to deal with non-compactness is to study manifolds with singularities. In 2014 Guillemin, Miranda and Pires \cite{Guillemin2014} introduced manifolds with a designated singular hypersurface, for which they coined the term \emph{$b$-manifold}. By restricting themselves to vector fields tangent to the singular hypersurface they were able to enlarge the space of differential forms to allow certain kinds of singularities. This makes it possible to generalize the symplectic and contact forms to \emph{$b$-symplectic} \cite{Guillemin2014} and \emph{$b$-contact} forms \cite{miranda2018,braddell2019}, respectively. These structures can be used to analyze Hamiltonian systems with non-compact energy hypersurfaces, where the singular hypersurface represents the behaviour of the system in 'the limit'. These structures have been shown to naturally arise in physical systems, for example in the analysis of the restricted three-body problem after the McGehee change of coordinates \cite{Delshams2019,Kiesenhofer2017,miranda2018,miranda2021}, when regularizing the collision of two particles in a certain non-gravitational potential \cite{braddell2019}, or in fluid dynamics \cite{Cardona2019, miranda2022}. Recently, Miranda and Oms \cite{miranda2021} formulated and analyzed the singular Weinstein conjecture and in the process found new periodic orbits on the non-compact energy level sets in the restricted $3$-body problem. Furthermore, Brugués, Miranda and Oms \cite{brugues2022} succeeded in defining Floer homology for some $b$-symplectic manifolds and used it to prove the Arnold conjecture in this case.

As one can see, in recent years there have been independently developed a variety of techniques to deal with the issue of non-compactness of a contact manifold. That rises obvious questions: do the techniques describe different objects? Are the classes of non-compact hypersurfaces, to which the different techniques apply, strictly disjoint? Are any of those techniques related? Do there exist objects to which those techniques can be applied alternatively? In this paper we will focus on two of the settings mentioned above namely the tentacular hyperboloids and the $b$-contact hypersurfaces. We show that in that case the last question has an affirmative answer. More precisely we show examples of tentacular hyperboloids that can be equipped with a $b$-contact structure.

Let us first introduce the class of Hamiltonians we will be working with.
Let $(M, \omega)$ be a symplectic manifold with a global Liouville vector field $Y$ and a Hamiltonian $H_0:M \to \R$, such that
\begin{equation}\label{defH0}
dH_0(Y)(x)-H_0(x)>0 \qquad \forall\ x \in M.
\end{equation}
We will consider Hamiltonians on $M \times \R^{2k}$ defined as
\begin{equation}\label{defH}
H(x,y):= H_0(x)+ y^T Ay,
\end{equation}
where $A$ is a symmetric, \emph{$J$-hyperbolic} matrix meaning that the matrix $\mathbb{J}A$ has no imaginary eigenvalues.

To formulate the main theorem we will first introduce some new nomenclature. We say that a $b$-symplectic manifold $(X, Z, \omega_b)$ is symplectically $b$-compatible with a symplectic manifold $(M, \omega)$ if every connected component of $X \setminus Z$ is symplectomorphic to $M$. For a given subset $A \subseteq M$ we can consider its images in $X \setminus Z$ through the symplectomorphisms between $M$ and the connected components of $X \setminus Z$. The closure in $X$ of the sum of those images will be called the $b$-image of $A$. Now we are ready to formulate the main result of this paper in the following theorem:

\begin{theorem}\label{thm_main}
Let $S \subseteq M \times \R^{2k}$ be a regular level set of a Hamiltonian $H$ as in \eqref{defH}. Then there exists a $b^3$-symplectic manifold $(X, Z, \omega_b)$, such that the $b$-image of $S$ in $X$ is a $b$-contact manifold. 
\end{theorem}
 
The following proposition shows that this theorem can be applied to a large class of tentacular hyperboloids (in fact to all known tentacular hyperboloids).
 
\begin{proposition}\label{prop:tent}
Let $A_0$ be a symmetric, positive definite, $2m$-by-$2m$ matrix and let $A_1$ be a symmetric, $J$-holomorphic, $2k$-by-$2k$ matrix. Then the Hamiltonian 
\begin{equation}\label{Hquad}
H(x,y):= \frac{1}{2}x^TA_0x+ \frac{1}{2}y^T A_1 y -1,
\end{equation}
satisfies \eqref{defH} and its zero level set is a tentacular hyperboloid.
\end{proposition}

To show that Hamiltonian $H_0(x):=\frac{1}{2}x^TA_0x-1$ satisfies \eqref{defH0} it suffices to take the global Liouville vector field $Y:=\frac{1}{2}x\partial_x$ and see that $dH_0(Y)-H_0\equiv 1$. Showing that $H$ is tentacular is non-trivial and is proven in Section \ref{sec_tent}.

The significance of the last result in the context of the Weinstein Conjecture for non-compact hypersurfaces will be more understandable if we provide a little more background on the class of the tentacular Hamiltonians and the application of Floer techniques to thereof. 
Rabinowitz Floer Homology is a symplectic invariant defined by Cieliebak and Frauenfelder in \cite{CieliebakFrauenfelder2009} associated to exact contact hypersurfaces in an exact convex symplectic manifold.
Its generators are the critical points of the so called Rabinowitz action functional for a suitable choice of defining Hamiltonian. This action functional was originally used by Rabinowitz in \cite{Rabinowitz1978} to prove the existence of periodic orbits on star-shaped hypersurfaces in $\R^{2n}$ and it differs from Floer’s original symplectic action functional by a Lagrangian multiplier, which forces the periodic solutions to lie on a prescribed energy surface. In fact, Rabinowitz Floer Homology is invariant under compact perturbations and is isomorphic to the singular homology if the hypersurface does not carry any contractible periodic orbits \cite{CieliebakFrauenfelder2009, pasquotto2018}. Those properties make Rabinowitz Floer homology a very suitable tool to study the Weinstein Conjecture.

However, before we make use of the properties of Rabinowitz Floer homology to answer the Weinstein Conjecture, we first have to calculate it.
The first computation of Rabinowitz Floer homology in the non-compact case was done by the second author together with Fauck and Merry in \cite{Fauck2021}. In Proposition \ref{prop:tent} we extend the class of tentacular hyperboloids. It turns out the result from \cite{Fauck2021} can also be applied to this larger class of tentacular hyperboloids defined in \eqref{Hquad} and can therefor be used to answer the Weinstein Conjecture for compact perturbations thereof.

The following Corollary is a straightforward application of \cite[Thm. 1.1]{Fauck2021} to the class of tentacular hyperboloids from Proposition \ref{prop:tent}:

\begin{corollary}
Let $S:=H^{-1}(0)$ be a tentacular hyperboloid defined as a zero level set of a tentacular Hamiltonian $H$ as in \ref{Hquad}. Then its Rabinowitz Floer homology is given by
$$
RFH_*(S):=\begin{cases}
\mathbb{Z}_2 & *= 1-m-k, -k,\\
0 & otherwise.
\end{cases}
$$
In particular, for $m>1$ we have $RFH_*(S) \neq H_{*+m+k-1}(S)$.
\end{corollary}

The Weinstein Conjecture is obvious for the tentacular hyperboloids described in Proposition \ref{prop:tent}. One of the strengths of Floer theoretical calculations, however, is that they are invariant under controlled perturbations. Therefore as an immediate consequence of the corollary above, we have:

\begin{corollary}
If $H$ is a Hamiltonian as in \eqref{Hquad} and $\{S_{\sigma}\}_{\sigma\in [0,1]}$ is a smooth $1$-parameter family of compact perturbations of $S:=H^{-1}(0)$ through strongly tentacular hypersurfaces, then each $S_{\sigma}$ carries a closed characteristics.
\end{corollary}

Theorem \ref{thm_main} together with \ref{prop:tent} bridges the two worlds of tentacular Hamiltonians and $b$-contact geometry. From \cite{pasquotto2018} we know that the Rabinowitz Floer homology for tentacular hyperboloids as in Proposition \ref{prop:tent} is well defined. In particular from \cite{pasquotto2017} we know that the Floer trajectories associated to the Rabinowitz action functional of a quadratic Hamiltonian as in \ref{Hquad} are contained in a compact set. On the other hand, Theorem \ref{thm_main} provides a $b$-symplectic manifold $(X, Z, \omega_b)$ and a pair of symplectomorphisms between the connected components of $(X\setminus Z)$ and $T^*\R^{m+k}$. Now consider the images of the Floer trajectories associated to the Rabinowitz action functional of $H$ via those symplectomorphisms. The images will be contained in a compact subset of $X$ and bounded away from $Z$. Consequently, one could use this result to generalize the Floer techniques to the $b$-symplectic setting and construct a Rabinowitz Floer homology for the $b$-contact type images of tentacular hyperboloids in the $b$-symplectic manifolds.

\subsection*{Organization of the paper}
In Section \ref{sec:prel} we recall the definition of tentacular Hamiltonians, symplectic hyperboloids and tentacular hyperboloids from \cite{pasquotto2018} and provide a brief introduction to $b$-symplectic geometry. In Section \ref{sec:wein} we provide examples of Hamiltonians satisfying condition \eqref{defH0} in arbitrary Weinstein manifolds. In Section \ref{sec_tent} we recall the H\"{o}rmander classification of quadratic forms under the action of the symplectic group $\operatorname{Sp}(\R^{2n})$, analyze different representatives of these equivalence classes, and prove Proposition \ref{prop:tent}. In Section \ref{sec_McG} we generalize the McGehee coordinate change to arbitrary dimension and then use it to construct a $b$-symplectic manifold which is symplectically $b$-compatible with the standard symplectic space $(\R^{2n}, \omega_0)$. Finally, in Section \ref{sec:bcont} we prove Theorem \ref{thm_main}. Additionally, in Appendix A we answer the questions of uniqueness of symmplectically $b$-compatible manifolds and in Appendix B we analyze the $b$-images of hypersurfaces.
\subsection*{Acknowledgments}

We would like to thank Eva Miranda for introducing us to the world of $b$-geometry through her many talks and for all the encouragement and interest shown in the subject of this paper.
We would also like to thank Cedric Oms for a fruitful discussion on the subject.
The first author would like to thank Barney Bramham for a helpful discussion about adjunction spaces.
The second author would like to thank Will Merry for all his support and mentoring during her stay at ETH Z\"{u}rich.

The first author is supported by the DFG under the Collaborative Research Center SFB/TRR 191 - 281071066 (Symplectic Structures in Geometry, Algebra and Dynamics).
The second author was partially supported by the SNF grant 200021\_182564\linebreak \href{http://p3.snf.ch/project-182564}{Periodic orbits on non-compact hypersurfaces}.
\setcounter{tocdepth}{1}
\tableofcontents
\section{Preliminaries}\label{sec:prel}

\subsection*{Tentacular Hamiltonians}The aim of this paper is to build a bridge between the worlds of tentacular Hamiltonians and $b$-contact geometry. Let us start with exploring the tentacular world first and recall the definition of tentacular Hamiltonians from \cite{pasquotto2018}.
%introduced by the second author in her PhD thesis \cite{Wisniewska2017}.

A continuous function $\phi:M \to \R$ is called \emph{coercive} if all its sublevel sets $\phi^{-1}((-\infty, a])$, $a \in \R$ are compact\footnote{Such a function is also sometimes called \emph{exhausting} by some authors.}. The \emph{Poisson bracket} of two functions $F$ and $H$ on a symplectic manifold $(M,\omega)$ is defined as $\{F,H\}:=\omega(X_F,X_H)$, where $X_F, X_H$ are Hamiltonian vector fields associated to $H$ and $F$ respectively, i.e. $dH = \omega(\ \cdot\ , X_H )$ and $dF = \omega( · , X_F )$. Recall that a \emph{Liouville vector field} is a vector field $Y$ on a symplectic manifold $(M, \omega)$, such that $d\iota_Y\omega=\omega$. A hypersurface $S\subseteq M$ is of \emph{contact type} if there exists a Liouville vector field transverse to $S$. A Liouville vector field $Y$ on the standard symplectic space $(\R^{2n},\omega_0)$ is called \emph{asymptotically regular} if $\|D Y(x)\| < c$ for some positive constant $c>0$ and all $x \in \R^{2n}$.

\begin{definition}\label{defTent}
Denote by $\mathcal{H}$ a class of smooth Hamiltonians $H : \R^{2n} \to \R$ satisfying the following axioms:
\begin{itemize}
\item[(h1)] there exist an asymptotically regular Liouville vector field $Y^\dagger$ and constants $c,c'>0$,
such that $dH(Y)(x) \ge c|x|^2 - c'$, for all $x\in \R^{2n}$;
\item[(h2)] (sub-quadratic growth) $\sup_{x\in \R^{2n}} \Vert D^3 H(x)\Vert \cdot |x| <\infty$;
\item[(h3)] there exists a coercive function $F$, defined in a neighborhood of $\Sigma$, such that for all $x\in \Sigma$ outside of a compact set either $\{H,F\}(x)\neq 0$ or $\{H,\{H,F\}\}(x)>0$. 
\end{itemize}\end{definition}

\begin{definition}
We say that $H^{-1}(0)$ is of \emph{restricted contact type} if there exists an asymptotically regular Liouville vector field $Y$, such
that $dH(Y )(x) > 0$ for all $x \in H^{-1}(0)$. Note that this implies that $H^{-1}(0)$ is a smooth hypersurface. We say that a Hamiltonian $H \in \mathcal{H}$ is \emph{strongly tentacular} if it's zero level set is of restricted contact type. If we drop the requirement that $F$ has to be coercive in axiom (h3), then we call $H$ simply \emph{tentacular}.
\end{definition}

\begin{remark}\label{rem:SpTent}
Note that $\mathcal{H}$ is closed under compact perturbations. In other words, if $H \in \mathcal{H}$ and $h\in C_c^\infty(\R^{2n})$, then $H+h \in \mathcal{H}$. Additionally, $\mathcal{H}$ is invariant under the action of the linear symplectic group $\operatorname{Sp}(\R^{2n})$. In other words if $H \in \mathcal{H}$ and $L \in \operatorname{Sp}(\R^{2n})$, then $H \circ L \in \mathcal{H}$.
\end{remark}

The tentacular Hamiltonians have been first studied in \cite{pasquotto2017} and \cite{pasquotto2018} where the second author together with Pasquotto and Vandervorst generalized the Rabinowitz Floer homology to this class.

\subsection*{Tentacular hyperboloids}
Let $A$ be a non-degenerate, quadratic, symmetric matrix and consider the non-degenerate quadratic Hamiltonian $H$ on $(\R^{2n},\omega_0)$
\begin{equation}
H(x)\coloneqq \frac{1}{2}x^TAx -1.\label{quadHam}
\end{equation}
The hypersurface $S\coloneqq H^{-1}(0)$ is diffeomorphic to $\s^{l-1}\times \R^{2n-l}$, where $(l,2n-l)$ is the signature of $A$. A hyperboloid is the zero level set of a quadratic Hamiltonian $H$ as in \eqref{quadHam} with $1\leq l  \leq 2n-1$.
Note that every hyperboloid is a hypersurface of contact type, as the radial Liouville vector field $Y=\frac{1}{2}x\partial_x$ satisfies $dH_x(Y)=H(x)+1=1\big|_S>0$.

In \cite[Sec. 9]{pasquotto2018} the second author together with Pasquotto and Vandervorst introduced the notion of \emph{symplectic hyperboloids}, which are equivalence classes of the set of hyperboloids under the action of the linear symplectic group $\operatorname{Sp}(\R^{2n})$. A symplectic hyperboloid is called a \emph{tentacular hyperboloid} if it admits a strongly tentacular representative of the form \eqref{quadHam}.

In \cite[Thm. 1.1]{Fauck2021} the second author together with Fauck and Merry computed the Rabinowizt Floer homology of a subset of tentacular hyperboloids and used it to answer the Weinstein Conjecture for quadratic perturbations thereof. In Proposition \ref{prop:tent} we enlarge the class of tentacular hyperboloids to all symplectic hyperboloids admitting a representative of the form \eqref{Hquad}, which enables us to apply \cite[Thm. 1.1]{Fauck2021} to this broader class.
\subsection*{\texorpdfstring{$b$}{b}-manifolds}
We will start by explaining the notion of a \emph{$b$-manifold}. We will follow the definitions introduced by Guillemin, Miranda and Pires in \cite{Guillemin2014}. These definitions were adapted from similar definitions introduced by Melrose in \cite{Melrose}. For a good introduction into this topic see also \cite{braddell2019}.
\begin{definition}
A couple $(M,Z)$ is called a \emph{$b$-manifold} if $M$ is an oriented smooth manifold and $Z$ an oriented hypersurface of $M$, which is called the \emph{singular hypersurface}. A vector field which is tangent to the singular hypersurface is called a \emph{$b$-vector field}. A \emph{$b$-map} is a map of $b$-manifolds $f:(M_1,Z_1)\to(M_2,Z_2)$ such that $f^{-1}(Z_2)=Z_1$ and $f$ is transverse to $Z_2$.
\end{definition}
Note that the set of $b$-vector fields is a strict subset of the space of vector fields. Since our end goal is to construct a generalization of symplectic and contact manifolds, it might seem counterproductive to think about a restriction of the space of vector fields. But in fact, this allows us to slightly enlarge the space of differential forms down the line.

On a (connected) $b$-manifold $(M,Z)$ we can find local coordinates $(x_1,\dots,x_n)$, such that $Z=x_1^{-1}(0)$. In these coordinates we have that the $b$-vector fields are (locally) generated by $\left(x_1\frac{\partial}{\partial x_1},\frac{\partial}{\partial x_2},\dots,\frac{\partial}{\partial x_n}\right)$. Thus the $b$-vector fields (locally) are a finitely generated free module. By the Serre-Swan-Theorem \cite[Thm. 11.32]{Jet}\ there exists a vector bundle, such that its sections are exactly the $b$-vector fields.
We call this vector bundle the \emph{$b$-tangent bundle} of $M$ and denote it by $\bb{TM}$. Its dual is the \emph{$b$-cotangent bundle} and is denoted $\btm$.

On a $b$-manifold $(M,Z)$ a \emph{$b$-differential form of degree k} is a section of $\bigwedge^k\left(\btm\right)$. These are denoted by $\bb{\Omega}^k(M)$.
Observe that if $Z\neq\emptyset$, then the $b$-differential forms contain the smooth differential forms as a strict subset. Note that the differential $d$ extends in a natural way to this larger space, hence we can also talk about closed and exact $b$-forms.

We can also apply the Serre-Swan theorem to the module which is locally generated by the vector fields $\left(x_1^m\frac{\partial}{\partial x_1},\frac{\partial}{\partial x_2},\dots,\frac{\partial}{\partial x_n}\right)$ for any natural number $m$. Then this gives us the \emph{$b^m$-tangent bundle} denoted $\bm{TM}$ and from this we can define as above its dual the \emph{$b^m$-cotangent bundle} denoted $\bmtm$ and lastly the \emph{$b^m$-differential forms} $\bm{\Omega}(M)$.
\begin{remark}\label{rem_jet}
Note that for $m>1$ the definition of the $b^m$-tangent bundle does in fact depend on the choice of the local coordinates. This can be resolved by specifying some jet-data of a defining function of the singular hypersurface, for a good reference to how this is done in detail see \cite{scott2016}. Fortunately, this will not be necessary for our purposes since we always specify local coordinates around the hypersurfaces for all our examples and thus implicitly use the corresponding jet-data.
\end{remark}
\subsection*{\texorpdfstring{$b$}{b}-symplectic manifolds}
Next we can define the $b$-symplectic manifolds the natural way following \cite[Def. 8]{Guillemin2014}:
\begin{definition}
A triple $(M,Z,\omega)$ is called a \emph{$b$-symplectic manifold}, if $(M,Z)$ is\linebreak a $b$-manifold and $\omega$ a closed $b$-2-form on $(M,Z)$ which is non-degenerate, i.e. its top power is non-vanishing as an element of the top $b$-forms.
\end{definition}
This definition (as most definitions from now on) could be extended the usual way to also define $b^k$-symplectic manifolds by using $b^k$-2-forms instead, but in the future this will mostly not be done explicitly to avoid clutter.
\begin{definition}\label{def_blvf}
On a $b$-symplectic manifold with $b$-symplectic form $\omega$, a $b$-vector field $X$ is a \emph{Liouville $b$-vector field}, if it satisfies $\mathcal{L}_X\omega=\omega$.
\end{definition}
\subsection*{\texorpdfstring{$b$}{b}-contact manifolds}
As in the smooth case, $b$-symplectic manifolds are always of even dimension, but there is also an odd-dimensional analogue, which were originally defined by Miranda and Oms in \cite[Def. 4.1]{miranda2018} the following way:
\begin{definition}
Let $(M, Z)$ be a $(2n+1)$-dimensional $b$-manifold. A \emph{$b$-contact structure} is the distribution given by the kernel of a $b$-form of degree 1, $\xi=\ker \alpha\subseteq \btm, \alpha \in \bb{\Omega}^1(M)$, that satisfies $\alpha\wedge(d\alpha)^n \neq 0$ as a element of $\bb{\Omega}^{2n+1}(M)$. We say that $\alpha$ is a \emph{$b$-contact form} and the pair $(M,\xi)$ a \emph{$b$-contact manifold}.
%A $(2n\!+\!1)$-dimensional $b$-manifold $(M,Z)$ together with a distribution is called a \emph{$b$-contact manifold}, if the distribution is (locally) given as the kernel of a $b$-1-form $\alpha$ and $\alpha\wedge(d\alpha)^n$ is non-vanishing (as a $b$-$(2n\!+\!1)$-form). The distribution is called the \emph{$b$-contact structure}.
\end{definition}
In this paper we will only talk about a special case of $b$-contact manifolds which arise as certain hypersurfaces of $b$-symplectic manifolds.
\begin{definition}\label{def_bchs}
Let $(M^{2n},Z,\omega)$ be a $b$-symplectic manifold. A hypersurface $S$ in $M$ is of \emph{$b$-contact type}, if it is transverse to Z and there exists a $b$-1-form $\alpha$ on $S$, such that $d\alpha=\omega\big|_S$ and $\alpha\wedge(d\alpha)^{n-1}$ is non-vanishing.
\end{definition}
If we consider a submanifold $N\subseteq M$ of a $b$-manifold $(M,Z)$, which is transverse to $Z$, then $(N,N\cap Z)$ is a $b$-manifold and there are natural inclusions $\bb{TN}\hookrightarrow\bb{TM}$ and $\bb{T^*N}\hookrightarrow\btm$. Consequently the transversality condition ensures that a $b$-contact type hypersurface is a $b$-manifold. In particular, a $b$-contact type hypersurface is a $b$-contact manifold. In general it is not obvious to determine which hypersurfaces are of $b$-contact type, but as in the smooth case we can characterize them using Liouville $b$-vector fields.
\begin{lemma}\label{lem_tlvf}
A hypersurface $S$ of a $b$-symplectic manifold is of $b$-contact type if and only if in the neighbourhood of $S$ there exists a Liouville $b$-vector field transverse to the hypersurface.
\end{lemma}
\begin{proof}
Let $S\subset M$ be a hypersurface of a $b$-symplectic manifold $(M,Z,\omega)$.
Let us assume that $V$ is a Liouville $b$-vector field which is defined locally around $S$ and is transverse to $S$. By assumption $V$ is tangent to $Z$ and transverse to $S$, hence $S$ has to be transverse to $Z$. Next we can define a 1-form on $S$ by setting $\alpha:=(\iota_V\omega)|_S$. Since $V$ is Liouville it satisfies
$$
d(\iota_V\omega)=\omega\qquad \text{and}\qquad \alpha\wedge(d\alpha)^{n-1}=\frac{1}{n}\iota_V \omega^n |_S \neq 0.
$$
Then $\alpha$ satisfies the requirements of Definition \ref{def_bchs} and thus makes $S$ into a $b$-contact type hypersurface. The proof in the other direction is done similarly as in the smooth case \cite[Lem. 2]{weinstein1979}.

\end{proof}

\subsection*{\texorpdfstring{$b$}{b}-compatible manifolds}
In this section we will introduce the concept of a \emph{$b$-compatible manifold}. The general idea is to consider a non-compact symplectic manifold and construct a $b$-symplectic manifold in such a way that the singular set of the $b$-symplectic manifold corresponds to sort of a 'compactification' of the original symplectic manifold\footnote{In fact the manifold we will end up with is not going to be compact, however the radial direction will be 'flipped' with $0$. For details see Section \ref{sec_bcm}.}. We hope that this construction will allow us to analyze the behavior of the symplectic manifold at infinity.
\begin{definition}
We say that a $b$-manifold $(X,Z)$ is \emph{$b$-compatible} with a smooth manifold $M$ if every connected component of $X\setminus Z$ is diffeomorphic to $M$. We say that a $b$-symplectic manifold $(X, Z, \omega_b)$ is \emph{symplectically $b$-compatible} with a symplectic manifold $(M, \omega)$ if every connected component of $X \setminus Z$ is symplectomorphic to $M$.
\end{definition}
\begin{remark}\label{rem_prod}
The concept of $b$-compatibility behaves well with products: Let $(M_1,\omega_1)$ and $(M_2, \omega_2)$ be two symplectic manifolds. If $(X,Z, \omega_b)$ is a $b$-symplectic manifold which is symplectically $b$-compatible to $(M_1,\omega_1)$ then $(X \times M_2, Z \times M_2, \omega_b \oplus \omega_2)$ is a $b$-symplectic manifold symplectically $b$-compatible with $(M_1 \times M_2, \omega_1 \oplus \omega_2)$.
\end{remark}

\begin{definition}\label{def_img}
Let $(X,Z)$ be a $b$-manifold $b$-compatible with a smooth manifold $M$. Let $\{C_i\}_{i\in I}$ be the connected components of $X \setminus Z$ and let $\varphi_i: M \to C_i$ be the associated diffeomorphisms. Then for every subset $S\subseteq M$ we define its \emph{$b$-image} in $X$ to be the set
$$
S_b:= \overline{\bigcup_{i\in I}\varphi_i(S)}.
$$
In other words, $S_b$ is the completion of the set $\bigcup_{i\in I}\varphi_i(S)$ in $X$.
\end{definition}
\begin{remark}
Note that the $b$-image heavily depends on the associated diffeomorphisms $\varphi_i$. Furthermore, since we can precompose any $\varphi_i$ by an arbitrary diffeomorphisms of $M$ to itself, there is in general no canonical choice of the $\varphi_i$.
\end{remark}

Assuming that you have a closed subset in the beginning, then the completion in Definition \ref{def_img} will in fact only add points belonging to $Z$. Using the same notation we can formulate this in the following way:
\begin{remark}
Consider $S\subseteq M$ a closed subset, then its $b$-image restricted to each connected component of $X\setminus Z$ is the same as its image under the associated diffeomorphism, i.e.:
\begin{equation*}
S_b\cap C_j=\phi_j(S)\qquad\forall j\in I.
\end{equation*}
This follows directly from the fact that taking the closure in $C_j$ is the same thing as restricting the closure in $X$ to $C_j$. 
\end{remark}
\section{Examples in Weinstein manifolds}\label{sec:wein}
The goal of this section is to construct examples of Hamiltonians satisfying \eqref{defH0} in Weinstein manifolds. These examples show that Theorem \ref{thm_main} can be applied to a large class of hypersurfaces, not only the tentacular hyperboloids. However, before we construct the examples of hypersurfaces in Weinstein manifolds, we will first analyze some basic geometrical properties of the hypersurfaces given as zero level sets of the Hamiltonians we will be working with defined in \eqref{defH}. More precisely, we will show that they are of contact type.

\begin{lemma}\label{lem:tentcont}
Let $H$ be a Hamiltonian as in \eqref{defH}. Then the hypersurface $S:=H^{-1}(0)$ is of contact type.
\end{lemma}
\begin{proof}
Consider the Liouville vector field $Y+\frac{1}{2}y\partial_y$ on $M\times T^*\R^k$. Then by \eqref{defH0} and \eqref{defH} for every $(x,y)\in S$  we have
\begin{align*}
dH\left(Y+\frac{1}{2}y\partial_y\right)(x,y) & =dH_0(Y)(x) + \frac{1}{2}dH_1(y \partial_y)=dH_0(Y)(x)+H_1(y)\\
& = H(x,y) +(dH_0(Y)-H_0)(x) = (dH_0(Y)-H_0)(x) >0.
\end{align*}
Consequently, the Liouville vector field $Y+\frac{1}{2}y\partial_y$ is everywhere transverse to $S$.
\end{proof}
\begin{remark}
Naturally, the zero level set of $H_0$ is also a contact type hypersurface in $(M,\omega)$ by \eqref{defH0}.
\end{remark}

In order to construct examples of Hamiltonians in Weinstein manifolds satisfying \eqref{defH0}, let us start with recalling some definitions. We say that a smooth function $\phi$ is called a \emph{Morse} function if all its critical points are non-degenerate. A vector field $Y$ on a manifold $M$ is said to be \emph{gradient-like} for a $C^1(M)$ function $\phi$ if for every Riemannian metric $g$ on $M$ there exists $\delta > 0$, such that
$$
d\phi(Y ) \geq \delta\left(\|Y \|^2 + \|\nabla \phi\|^2\right).
$$
\begin{definition}
A \emph{Weinstein manifold} $(M, \omega, Y, \phi)$ is a symplectic manifold $(M, \omega)$ with a complete Liouville vector field $Y$, which is gradient-like for a coercive Morse function $\phi: M \to \R$. A Weinstein manifold $(M, \omega, Y, \phi)$ is said to be \emph{of finite type} if $\phi$ has finitely many critical points.
\end{definition}

\begin{lemma}\label{lem:H0}
Let $(M, \omega, Y, \phi)$ be a Weinstein manifold of finite type. Then there exists a Hamiltonian function $H_0: M \to \R$ satisfying \eqref{defH0}.
\end{lemma}
Before proving the lemma above we will need another one, which shows that a Weinstein manifold of finite type is a strong deformation retract of its compact subset.
\begin{lemma}\label{lem:tau}
Let $(M, \omega, Y, \phi)$ be a Weinstein manifold of finite type. Then there exists $a\in \R$, such that $\phi^{-1}((-\infty, a])$ is a strong deformation retract of the manifold $M$.
\end{lemma}

\begin{proof}
Let us denote by $\psi^t$ the flow of $Y$.
By assumption $\phi$ is coercive and has finitely many critical points, hence there exists $a\in \R$, such that $\operatorname{Crit}\phi \subseteq \phi^{-1}((-\infty, a))$.  
We would like to show that there exists a smooth function $\tau: \phi^{-1}([a,+\infty)) \to [0,+\infty)$, such that $ \phi \circ \psi^{-\tau(x)}(x)=a$ for all $x\in \phi^{-1}([a,+\infty))$.

Let us fix $x_0\in  \phi^{-1}([a,+\infty))$ and denote $a':=\phi(x_0)$. By assumption $\phi$ is a coercive function, hence $\phi^{-1}([a,a'])$ is a compact set. On the other hand, $Y$ is gradient-like for $\phi$, hence $d\phi(Y)>0$ on $\phi^{-1}([a,\infty))$. Denote
$$
\delta:= \inf \left\lbrace d\phi(Y)(x)\ \big|\ a'\leq \phi(x)\leq a\right\rbrace>0.
$$
Consequently, for all $t \in \R_+$, such that $0 \leq \phi(\psi^{-t}(x_0))$ we have
$$
\phi\circ \psi^{-t}(x_0) = a' - \int_{-t}^0 dH(Y)(\psi^{s}(x_0))ds \leq a'-\delta t.
$$
In other words, $\phi$ linearly decreases along the flow line of $-Y$. Consequently, there exists $t_0\in \R_+$, such that $\phi \circ \psi^{-t_0}(x_0)=a$. This assures that the function 
\begin{align}
\tau: \phi^{-1}([a,+\infty)) &\to [0,+\infty),\nonumber\\
\textrm{such that} \qquad \phi \circ \psi^{-\tau(x)}(x)& =a\qquad \forall\ x\in  \phi^{-1}([a,+\infty)),\label{tau}
\end{align}
 is well defined. Note that $\tau\big|_{\phi^{-1}(a)} \equiv 0$.
 
Consequently, the map $r_0: M\times [0,1] \to M$ defined 
$$
r_0(x,t):=\begin{cases}
\psi^{-t\tau(x)}(x) & \textrm{for}\quad x\in \phi^{-1}([a,+\infty)),\\
x &  \textrm{for}\quad x\in \phi^{-1}((-\infty, a]),\end{cases}
$$
is a strong deformation retract of $M$.
\end{proof}

Now we are ready to prove Lemma \ref{lem:H0}:
\vspace*{.25cm}

\noindent\textit{Proof of Lemma \ref{lem:H0}:}
Let $\psi^t$ denote the flow of the Liouville vector field. Since the Weinstein manifold is of finite type there exists a constant $a\in \R$, such that $\operatorname{Crit \phi} \subseteq \phi^{-1}((-\infty,a])$. By Lemma \ref{lem:tau} there exists a smooth function $\tau: \phi^{-1}((a,\infty)) \to \R_+$, such that
\begin{equation}\label{phi}
\phi \circ \psi^{-\tau(x)}(x)=a\quad \textrm{and} \quad \tau \circ \psi^t (x)=t+\tau(x)\qquad \forall\ x\in \phi^{-1}((a,\infty)).
\end{equation}
Now let $\chi: \R \to [0,1]$ be a smooth, non-decreasing function, such that $\chi(t)\equiv 0$ for $t\leq 0$, $\chi(t)\equiv 1$ for $t\geq 1$ and $\chi'(t)>0$ on $(0,1)$.
Denote $a':=\max\{ \phi \circ \psi^1(x)\ |\ x\in \phi^{-1}(a)\}$.
We define our Hamiltonian function
$$
H_0(x):=(\phi(x)-a')(1-\chi \circ \tau (x))+(\exp \circ \tau(x)-1) \cdot \chi \circ \tau (x).
$$
{Using \eqref{phi} we substitute $t:=\tau(x), y:=\psi^{-\tau(x)}(x)\in \phi^{-1}(a)$ for every $x\in \phi^{-1}([a,\infty))$ and calculate
\begin{align*}
H_0 \circ \psi^t (y) & = (\phi \circ \psi^t(y)-a')(1-\chi(t))+(\exp(t)-1)\chi(t),\\
\frac{d}{dt}\left( H_0 \circ \psi^t \right) & = d\phi (Y)(1-\chi(t))+\exp(t)\chi(t)+\chi'(t)(\exp(t)-1+a'-\phi\circ \psi^t(y)),\\
(dH_0(Y)-H_0)\circ \psi^t(y) & = (d\phi (Y)+a'-\phi\circ \psi^t(y))(1-\chi(t))+\chi(t)\\
& \hspace*{4.5cm}+\chi'(t)(\exp(t)-1+a'-\phi\circ \psi^t(y))\\
dH_0(Y)-H_0(x) & = (d\phi (Y)+a'-\phi(x))(1-\chi\circ \tau (x))+\chi\circ \tau (x)\\
& \hspace*{4.5cm}+\chi'\circ \tau (x)(\exp\circ \tau(x)-1+a'-\phi(x)).
\end{align*}
Using the last equation we will now show that $dH_0(Y)-H_0>0$ everywhere on $M$.

Recall that by construction $a'>\phi(x)$ for all $x \in \phi^{-1}((-\infty, a])\cup \tau^{-1}((0,1))$. On the other hand, by assumption $\phi$ is gradient-like for $Y$, hence $d\phi(Y)\geq 0$ on all of $M$. Consequently, $(d\phi (Y)+a'-\phi(x))(1-\chi\circ \tau(x))>0$ for all $x \in \phi^{-1}((-\infty, a])\cup \tau^{-1}((0,1))$ and zero otherwise.

By construction, both $\chi'(t)$ and $\exp(t)-1$ are positive on $(0,1)$. Hence the last summand $\chi'\circ \tau (x)(\exp\circ\tau(x)-1+a'-\phi(x))$ is positive for all $x\in\tau^{-1}((0,1))$ and also zero otherwise.

Finally, $\chi(t)$ is by construction positive for all $t>0$ and non-negative for all $t$. We can conclude that all summands are non-negative on the whole $M$ and at least two of the summands are positive on either of the sets $\phi^{-1}((-\infty, a])\cup \tau^{-1}((0,1))$ and $\tau^{-1}((0,+\infty))$ covering $M$. Therefore, $dH_0(Y)-H_0>0$ everywhere on $M$.

\hfill $\square$
\section{Symplectic hyperboloids}\label{sec_tent}

The main aim of this section will be to prove Proposition \ref{prop:tent}, i.e. show that Hamiltonians of the form \eqref{Hquad} are tentacular in the sense of Definition \ref{defTent}. Recall that if a Hamiltonian $H$ is tentacular, then for all linear symplectomorphisms $L \in \operatorname{Sp}(\R^{2n})$ the Hamiltonian $H \circ L$ is also tentacular. Therefore, prior to proving Proposition \ref{prop:tent}, we will first discuss different representatives of the equivalence classes of the set of quadratic Hamiltonians under the action of the linear symplectic group $\operatorname{Sp}(\R^{2n})$ in search of a representative for which proving the Axioms of Definition \ref{defTent} would be the easiest.

\subsection{H{\"o}rmander classification}

First, we will recall the H{\"o}rmander classification of quadratic forms \cite[Thm.\ 3.1]{Hormander1995}. Let $A$ be a $2n$-by-$2n$, non-degenerate, symmetric matrix. Observe that if $\lambda$ is an eigenvalue of $\mathbb{J}A$ with a corresponding $m\times m$ block in the Jordan decomposition, then its additive inverse and its complex conjugate are also eigenvalues of $\mathbb{J}A$ each with a $m\times m$ block.

By the H{\"o}rmander classification, there exists a symplectic change of coordinates, after which we can assume that $T^*\R^n$ splits into a direct sum of symplectic subspaces $S_i$,  such that the matrix $A$ is a diagonal block matrix consisting of matrices $A_i=A|_{S_i}$. Each matrix $A_i$ corresponds to a $m\times m$ block with an eigenvalue $\lambda$ in the Jordan decomposition of $\mathbb{J}A$ and it is determined as follows:
\begin{enumerate}[label=(\alph*)]
\item if $\Im(\lambda)=0$, then $S_i\coloneqq T^*\R^{m}$ and $A_i\coloneqq\left(\begin{smallmatrix} 0 & B\\ B^T & 0\end{smallmatrix}\right)$ is a $2m \times 2m$ block matrix where 
\begin{equation}\label{eq_Ba}
B:=\left(\begin{array}{c c c c}
\mu\\
1 & \mu\\
& \ddots & \ddots \\
 & & 1 & \mu
\end{array}\right),
\end{equation}
with $\mu:=|\Re (\lambda)|$.
In other words, $B$ is a $m \times m$ matrix with $\mu$'s on the diagonal and with $1$'s under the diagonal for $m>1$. 
The signature of $A_i$ is $(m,m)$ and the corresponding Hamiltonian $H_i \colon T^*\R^{m} \to \R$ is given by
\[
H_i(q,p)\coloneqq \mu\sum_{j=1}^{m}q_jp_j+\sum_{j=1}^{m-1}q_{j+1}p_j.
\]
\item if $\Re(\lambda)\neq 0,\ \Im(\lambda)\neq 0$, then $S_i\coloneqq T^*\R^{2m}$ and
 $A_i\coloneqq  \left(\begin{smallmatrix} 0 & B\\ B^T & 0\end{smallmatrix}\right)$ is an $4m\times 4m$ block matrix where
\begin{equation}\label{eq_Bb}
B\coloneqq \left(\begin{array}{c c c c c c c }
\mu_1 & \mu_2 & 1 \\
-\mu_2 & \mu_1  & & 1\\
&  & \mu_1 & \mu_2 & 1 \\
 & & - \mu_2 & \mu_1 &  & \ddots\\
 & & &  & \ddots &  & 1 \\
 & & & &  &\mu_1 & \mu_2 \\
 & & & &  & -\mu_2 & \mu_1 \\
\end{array}\right),
\end{equation}
with $\mu_1\coloneqq |\Re(\lambda)|, \mu_2\coloneqq |\Im(\lambda)|$.
The signature of $A_i$ is $(2m,2m)$ and the Hamiltonian $H_i \colon T^*\R^{2m} \to \R$ is given by
\[
H_i(q,p)\coloneqq   \sum_{j=1}^{2m-2}q_jp_{j+2}+\mu_1\sum_{j=1}^{2m}q_jp_j+\mu_2\sum_{j=1}^{m} \big( q_{2j}p_{2j-1}-q_{2j-1}p_{2j} \big).
\]

\item if $\Re(\lambda)=0$, then $S_i\coloneqq T^*\R^m$ and $A_i\coloneqq \left(\begin{smallmatrix} B & 0\\ 0 & B^P\end{smallmatrix}\right)$ is an $2m\times 2m$ block matrix, where 
$$
B\coloneqq \gamma \left(\begin{array}{c c c c c}
& &&& \mu\\
& & & \mu & -1\\
&& \reflectbox{$\ddots$} & \reflectbox{$\ddots$}\\
& \mu & - 1\\
\mu & -1 \\
\end{array}\right)
$$
for $\mu\coloneqq |\Im(\lambda)|$, some $\gamma=\pm 1$ and $B^P$ being the reflection of $B$ with respect to the anti-diagonal. That is, if $B\coloneqq \{b_{j,k}\}_{j,k=1}^m$ then $B^P\coloneqq \{b_{\delta(j,k)}\}_{j,k=1}^m$ with\linebreak $\delta(j,k)=(m{+}1\,{-}j, m{+}1\,{-}k)$.
The signature of $A_i$ is $(m,m)$ if $m$ is even\linebreak and $(m+\gamma,m-\gamma)$ if $m$ is odd, and the corresponding Hamiltonian $H_i \colon T^*\R^m \to \R$ is equal to
\[
H_i(q,p)\coloneqq \frac{\gamma}{2}\left(\mu\sum_{j=1}^m \big(q_jq_{m+1-j}+ p_jp_{m+1-j} \big)-\sum_{j=1}^{m-1} \big(q_{j+1}q_{m+1-j} +p_jp_{m-j} \big)\right).
\]
\end{enumerate}

The H{\"o}rmander classification determines a unique (up to permutation of the blocks) representative of the equivalence class of non-degenerate matrices under the action of $\operatorname{Sp}(\R^{2n})$. In particular, after a linear symplectic change of coordinates, every quadratic Hamiltonian of the form  $H(x):=x^TAx-1$ can be represented as a sum of Hamiltonians $H_i \colon S_i\to\R,$
\[
H(x)=\sum H_i(x_i)-1, \qquad H_i(x_i)\coloneqq  \frac{1}{2} x_i^TA_ix_i, \qquad x_i\in S_i,
\]
with $A_i$ either of type (a), (b) or (c).

\subsection{$J$-hyperbolic matrices}
In this subsection we will discuss the class of $J$-hyperbolic matrices. Both the $J$-hyperbolic and the tentacular conditions are invariant under the action of the $\operatorname{Sp}(\R^{2n})$ group. However, the $J$-hyperbolic condition is easy to verify regardless of the choice of the representative, whereas to prove the tenatacular condition it would be more convenient to work with other representatives. Therefore, in this subsection we will investigate other representatives of the set of $J$-hyperbolic matrices under the action of the $\operatorname{Sp}(\R^{2n})$ group.

\begin{proposition}\label{prop_JAB}
Consider a $2k \times 2k$ symmetric, nondegenerate matrix $A$. Then the following conditions are equivalent:
\begin{enumerate}[label=\roman*)]
\item The matrix $A$ is $J$-hyperbolic.
\item There exists a $k\times k$ matrix $B$, such that $B+B^T$ is positive definite and $A$ is simplectically conjugate to the matrix
$$
\left(\begin{array}{c c}
0 & B\\
B^T & 0
\end{array} \right).
$$
\end{enumerate}
\end{proposition}

Unfortunately, we cannot prove Proposition \ref{prop_JAB} by simply applying the H\"{o}rmander classification, because the matrices $B$ of type \eqref{eq_Ba} and \eqref{eq_Bb} do not always satisfy the condition that $B+B^T$ is positive definite. Therefore we will need to find more suitable representatives, which we will do in the following proposition:
\begin{proposition}\label{prop_Beps}
Let $A$ be a $2n \times 2n$, non-degenerate, symmetric matrix, such that $\mathbb{J}A$ has no imaginary eigenvalues. Then for every $\varepsilon>0$ there exists a symplectic change of coordinates, after which we can assume that $T^*\R^n$ splits into a direct sum of symplectic subspaces $S_i$,  such that the matrix $A$ is a diagonal block matrix consisting of matrices\linebreak $A_i=A|_{S_i}$. Each matrix $A_i$ corresponds to a $m\times m$ block with an eigenvalue $\lambda$ in the Jordan decomposition of $\mathbb{J}A$ and it is determined as follows:
\begin{enumerate}[label=(\alph*)]
\item if $\Im(\lambda)=0$, then $S_i\coloneqq T^*\R^{m}$ and $A_i\coloneqq\left(\begin{smallmatrix} 0 & B_\varepsilon\\ B_\varepsilon^T & 0\end{smallmatrix}\right)$ is a $2m \times 2m$ block matrix where 
\begin{equation}\label{eq_Ca}
B_\varepsilon \coloneqq \begin{pmatrix}\mu \\ \varepsilon & \mu \\ & \ddots & \ddots \\ && \varepsilon& \mu \end{pmatrix}
\end{equation}
with $\mu:=|\Re (\lambda)|$.
In other words, $B_\varepsilon$ is an $m \times m$ matrix with $\mu$'s on the diagonal and with $\varepsilon$'s under the diagonal for $m>1$. 
\item if $\Re(\lambda)\neq 0$ and $\Im(\lambda)\neq 0$, then $S_i\coloneqq T^*\R^{2m}$ and
 $A_i\coloneqq  \left(\begin{smallmatrix} 0 & B_\varepsilon\\ B_\varepsilon^T & 0\end{smallmatrix}\right)$ is an $4m\times 4m$ block matrix where
\begin{equation}\label{eq_Cb}
B_\varepsilon \coloneqq \begin{pmatrix}
\mu_1 & \mu_2 & \varepsilon\\ -\mu_2 & \mu_1 & & \varepsilon\\
&& \mu_1 & \mu_2 & \varepsilon \\
&& -\mu_2 & \mu_1 && \ddots\\
&&&& \ddots && \varepsilon\\
&&&&& \mu_1 & \mu_2\\
&&&&& -\mu_2 & \mu_1
\end{pmatrix}
\end{equation}
with $\mu_1\coloneqq |\Re(\lambda)|, \mu_2\coloneqq |\Im(\lambda)|$.
\end{enumerate}
\end{proposition}
The proof of the above proposition is obtained by combining the H\"{o}rmander classification result \cite[Thm.\ 3.1]{Hormander1995} and the two lemmas presented below:

\begin{lemma}
A matrix $B$ as in \eqref{eq_Ba} is similar to a matrix $B_\varepsilon$ as in \eqref{eq_Ca} for any $\varepsilon>0$, whereas a matrix $B$ as in \eqref{eq_Bb} is similar to a matrix $B_\varepsilon$ as in \eqref{eq_Cb} for any $\varepsilon>0$.
\end{lemma}

\begin{proof}
Let us first consider $B$ to be an $m\times m$ matrix as in \eqref{eq_Ba}. For $\varepsilon>0$ we define an $m\times m$ diagonal matrix $D_\varepsilon$ with numbers $\{\varepsilon^m, \varepsilon^{m-1}, \dots, 1\}$ on the diagonal. Naturally, its inverse $D_\varepsilon^{-1}$ will have  $\{\varepsilon^{-m}, \varepsilon^{1-m}, \dots, 1\}$ on the diagonal. Let us calculate
$$
D_\varepsilon^{-1}B D_\varepsilon \coloneqq 
\begin{pmatrix}\epsilon^{-m}\hspace{-.3cm}\\&\epsilon^{1-m}\hspace{-.9cm}\\&&\ddots\\&&&1\end{pmatrix}
\begin{pmatrix}\mu \\ 1 & \mu \\ & \ddots & \ddots \\ && 1& \mu \end{pmatrix}
\begin{pmatrix}\epsilon^m\hspace{-.3cm}\\&\epsilon^{m-1}\hspace{-.9cm}\\&&\ddots\\&&&1\end{pmatrix}
= \begin{pmatrix}\mu \\ \varepsilon & \mu \\ & \ddots & \ddots \\ && \varepsilon& \mu \end{pmatrix}=B_\varepsilon,
$$
which proves our claim.

Let us now consider $B$ to an $2m\times 2m$ matrix as in \eqref{eq_Bb}. For $\varepsilon>0$ we define a $2m\times 2m$ diagonal matrix $D_\varepsilon$ with numbers $\{\varepsilon^m, \varepsilon^m, \varepsilon^{m-1}, \varepsilon^{m-1}, \dots, 1, 1\}$ on the diagonal. Naturally, its inverse $D_\varepsilon^{-1}$ will have  $\{\varepsilon^{-m}, \varepsilon^{-m}, \varepsilon^{1-m}, \varepsilon^{1-m}, \dots,1, 1\}$ on the diagonal. Analogously as in the previous case we can verify that $B_\varepsilon=D_\varepsilon B D_\varepsilon^{-1}$ will be of the form as in \eqref{eq_Cb}.
\end{proof}

\begin{lemma}\label{lem_simpB}
Suppose that the two $m \times m$ matrices $B$ and $C$ are similar. Then the matrices
$$
\begin{pmatrix} & B\\ B^T \end{pmatrix} \qquad \text{and} \qquad \begin{pmatrix}& C\\ C^T\end{pmatrix},
$$
are simplectically conjugate.
\end{lemma}
\begin{proof}
First observe that for every $D\in \operatorname{GL}(\R^m)$ the matrix $\left(\begin{smallmatrix} \left(D^{-1}\right)^T & \\ & D \end{smallmatrix}\right)$
is symplectic. Indeed, we can calculate
$$
\begin{pmatrix}D^{-1}\\ & D^T\end{pmatrix}
\begin{pmatrix} & \operatorname{Id}\\ -\operatorname{Id}\end{pmatrix}
\begin{pmatrix}\left(D^{-1}\right)^T\\&D\end{pmatrix}=
\begin{pmatrix} & D^{-1}\\ -D^T\end{pmatrix}
\begin{pmatrix}\left(D^{-1}\right)^T\\&D\end{pmatrix}=
\begin{pmatrix} & \operatorname{Id}\\ -\operatorname{Id}\end{pmatrix}.
$$

By assumption $B$ and $C$ are similar, hence there exists a matrix $D\in \operatorname{GL}(\R^m)$, such that $C=D^{-1}BD$. What is left is to calculate:
\begin{align*}
\begin{pmatrix}D^{-1}\\&D^T\end{pmatrix}
\begin{pmatrix} & B\\ B^T\end{pmatrix}
\begin{pmatrix}\left(D^{-1}\right)^T\\ & D\end{pmatrix} & =
\begin{pmatrix}D^{-1}\\&D^T\end{pmatrix}
\begin{pmatrix} & BD\\ B^T \left(D^{-1}\right)^T\end{pmatrix}\\
= \begin{pmatrix}&D^{-1}BD\\ D^TB^T\left(D^{-1}\right)^T\end{pmatrix} &=
\begin{pmatrix} & C\\ C^T\end{pmatrix}
\end{align*}
\end{proof}

Having proven Proposition \ref{prop_Beps} we are finally ready to prove Proposition \ref{prop_JAB}:

\textit{Proof of Proposition \ref{prop_JAB}}:\\
i) $\Rightarrow$ ii):\\
By assumption the matrix $\mathbb{J}A$ has no imaginary eigenvalues, hence by Proposition \ref{prop_Beps} the matrix $A$ is symplectically conjugate to a matrix of the form $\left(\begin{smallmatrix} 0 & B\\ B^T & 0\end{smallmatrix}\right)$, where $B$ itself is a block matrix consisting of matrices $B_\varepsilon$ each corresponding to an eigenvalue $\lambda \in \sigma(\mathbb{J}A)$. In case $\Im(\lambda)=0$ then $B_\varepsilon$ is of the form as in \eqref{eq_Ca}, whereas in case $\Im(\lambda)\neq 0$ then $B_\varepsilon$ is of the form as in \eqref{eq_Cb}. Note that if $B_\varepsilon$ is of the form as in \eqref{eq_Ca} or of the form as in \eqref{eq_Cb}, then
$$
B_\varepsilon+B_\varepsilon^T = \begin{pmatrix} 2\mu & \varepsilon\\ \varepsilon & \ddots & \ddots\\ &\ddots & &\varepsilon\\ &&\varepsilon & 2\mu \end{pmatrix} \qquad \text{or} \qquad
B_\varepsilon+B_\varepsilon^T =  \begin{pmatrix} 2\mu_1 & & \varepsilon\\  & 2\mu_1 & & \ddots\\ \varepsilon &&  \ddots && \varepsilon\\ & \ddots& &2 \mu_1 & \\ && \varepsilon &  & 2\mu_1 \end{pmatrix},
$$
respectively. In particular, in both cases $B_\varepsilon+B_\varepsilon^T$ is diagonally dominant, whenever $\varepsilon < \mu, \mu_1$. Furthermore, whenever $\varepsilon < |\Re(\lambda)|$, then by the Levy-Desplanques theorem \cite{levy1881}, \cite{desplanques1887} the corresponding matrix $B_\varepsilon+B_\varepsilon^T$ is positive definite.

Consequently, $B+B^T$ will be a block matrix with matrices $B_\varepsilon+B_\varepsilon^T$ on the diagonal each corresponding to an eigenvalue $\lambda \in \sigma(\mathbb{J}A)$.  Now, if we take $\varepsilon < \min_{\lambda\in \sigma(\mathbb{J}A)}|\Re(\lambda)|$, then by the Levy-Desplanques theorem \cite{levy1881}, \cite{desplanques1887} we can conclude that $B+B^T$ itself is positive definite.

ii) $\Rightarrow$ i):\\
We will first prove that if $A\coloneqq \left(\begin{smallmatrix} 0 & B\\ B^T & 0\end{smallmatrix}\right)$ and $B+B^T$ is positive definite, then $\mathbb{J}A=\left(\begin{smallmatrix} B^T & 0\\ 0 & -B\end{smallmatrix}\right)$ has no imaginary eigenvalues. Observe that $\sigma(\mathbb{J}A)=\sigma(B^T)\cup \sigma (-B)=\sigma(B)\cup(- \sigma(B))$, so it suffices to prove that $B$ has no imaginary eigenvalues. We will prove it by contradiction.

Suppose that $B$ has an imaginary eigenvalue, i.e. there exists $\mu \in \R\setminus \{0\}$, such that $i\mu \in \sigma(B)$. In other words, there exists a complex vector $v$, such that $Bv=i\mu v$. Then the real vectors $v_1\coloneqq \Re(v)$ and $v_2\coloneqq \Im(v)$ satisfy
$$
Bv_1=-\mu v_2 \qquad \text{and}\qquad Bv_2= \mu v_1.
$$
From $v$ and $\mu$ being non-zero we can deduce that both $v_1$ and $v_2$ are non-zero. We can calculate:
\begin{align*}
\frac{1}{2} v_1^T(B+B^T)v_1 & = v_1^TBv_1=-\mu v_1^Tv_2,\\
\frac{1}{2} v_2^T(B+B^T)v_2 & = v_2^TBv_2 = \mu v_2^Tv_1.
\end{align*}
By assumption $B+B^T$ is positive definite, but the two numbers on the right hand side of the equation cannot be both positive at the same time, which brings us to a contradiction.

\hfill $\square$

\subsection{Tentacular hyperboloids}
Our aim now is to prove Proposition \ref{prop:tent}. In other words we want to prove that the zero level set of Hamiltonian \eqref{Hquad} is a tentacular hyperboloid. Naturally, all quadratic Hamiltonians satisfy Axiom (h2) of Definition \ref{defTent}. However, the other properties of tentacular Hamiltonians are slightly more challenging to verify. Therefore, we will use different representatives of the equivalence classes under the action of the linear symplectic group $\operatorname{Sp}(\R^{2n}$) to verify those conditions.

Before proving Proposition \ref{prop:tent} we will first prove some lemmas:

\begin{lemma}\label{lem:h1}
Let $A$ be an $2m$-by-$2m$ symmetric, positive definite matrix and let $B$ be a $k$-by-$k$ matrix, such that $B+B^T$ is positive definite. Define
$$
H(x, q, p):= \frac{1}{2}x^TAx + p^TBq-1\qquad x\in T^*\R^m,\ (q,p)\in T^*\R^k.
$$
Then $H$ satisfies Axiom (h1) of Definition \ref{defTent} and its zero level set is of contact type.
\end{lemma}
\begin{proof}
By assumption the matrices $A$ and $B+B^T$ are positive definite, so there exists a positive constant $c>0$, such that 
$$
x^TAx \geq c |x|^2 \quad \forall\ x\in T^*\R^m\qquad \textrm{and}\qquad v^T(B+B^T)v \geq c|y|^2 \quad \forall\ v\in \R^k.
$$
For a constant $a := \frac{\|B\|}{c}+1$ we define a Liouville vector field
$$
Y(x,q,p):= \frac{1}{2}(x\partial_x+q\partial_q+p\partial_p)+a(q\partial_p+p\partial_q).
$$
Now we can calculate
\begin{align*}
dH (x, q,p) & = x^T A dx + \frac{1}{2}\left( p^TB dq + q^TB^Tdp\right)\\
dH\left(Y\right) & = \frac{1}{2}x^TAx+p^TBq+a (q^TBq+p^TBp)\\
 & = \frac{1}{2}x^TAx+p^TBq+\frac{a}{2}\left(q^T(B+B^T)q+p^T(B+B^T)p\right)\\
 & \geq \frac{c}{2}|x|^2 - \|B\||q||p|+\frac{a c}{2}\left(|q|^2+|p|^2\right)\\
 & \geq \frac{1}{2}\left(c|x|^2 + \left(a c-\|B\|\right)\left(|q|^2+|p|^2\right)\right)\\
 & \geq \frac{c}{2} \left(|x|^2+|q|^2+|p|^2\right).
\end{align*}

This proves that $H$ satisfies Axiom (h1). Moreover, since $0 \notin H^{-1}(0)$ we also obtain that $H^{-1}(0)$ satisfies the contact type condition.
\end{proof}

\begin{lemma}\label{lem:h3}
Let $r_1, \dots r_m>0$ be a sequence of positive constants and let $B$ be a $k$-by-$k$ matrix, such that $B+B^T$ is positive definite. For $(x,y)=(x_1, \dots x_m, y_1, \dots y_m,)\in T^*\R^m$ and $(q,p)\in T^*\R^k$ define
$$
H(x,y, q, p)  := \sum_{i=1}^{m}r_i\left(x_i^2+y_i^2\right)+ p^TBq-1
$$
Then $H$ satisfies Axiom (h3) of Definition \ref{defTent}.
\end{lemma}
\begin{proof}
Let us define a coercive function
$$
F(x,y, q, p) := \frac{1}{2}\left(|x|^2+|y|^2+|q|^2 + |p|^2\right).
$$
Let us calculate $\{H,F\}=dF(X_H)$, where $X_H$ is the Hamiltonian vector field of $H$.
\begin{align*}
dF (x,y,q,p) & = xdx+ydy+qdq+pdp,\\
X_H(x,y,q,p) & = \sum_{i=1}^m r_i\left(y_i\partial_{x_i}-x_i\partial_{y_i}\right)+q^TB^T\partial_q-p^TB\partial_p,\\
\{H,F\}(x,y,q,p) &= dF(X_H) = q^TBq+p^TBp.
% d \{H,F\}(x,y,q,p) & = q^T(B+B^T)dq+p^T(B+B^T)dp,\\
%\{H, \{H, F\}\} (x,y,q,p) & = d\{H, F\}(X_H)= 
\end{align*}
Note that the first part of the Hamiltonian gives no input to $\{H,F\}$ - it disappears!

By assumption $B+B^T$ is positive definite, so there exists a positive constant\linebreak $0 < c \leq \min_i r_i$, such that $v^T(B+B^T)v \geq c|v|^2$ for all $v\in \R^k$.
Set $\delta:= \frac{c}{c+\|B\|}$ and estimate
\begin{align*}
\left(\{H,F\}+\delta H\right)(x,y,q,p) & = \frac{1}{2}\left(q^T(B+B^T)q+p^T(B+B^T)p\right)+\delta \sum_{i=1}^{m}r_i\left(x_i^2+y_i^2\right)+\delta p^TBq - \delta\\
& \geq \frac{c}{2}\left(\delta \left(|x|^2+|y|^2\right)+|q|^2+|p|^2\right)-\delta \|B\||p||q| -\delta\\
& \geq \frac{c \delta}{2}\left(|x|^2+|y|^2\right) +\frac{1}{2}(c- \delta\|B\|)\left(|q|^2+|p|^2\right)-\delta\\
& = \frac{c\delta}{2}\left(|x|^2+|y|^2+|q|^2+|p|^2\right)-\delta
\end{align*}
From the inequality above we can deduce that $\{H,F\}>0$ on $H^{-1}(0)$ outside a compact set $F^{-1}((-\infty, \frac{1}{c}])$. Thus $H$ satisfies Axiom (h3).
\end{proof}

\noindent\textit{Proof of Proposition \ref{prop:tent}:}
By Preposition \ref{prop_JAB} for every Hamiltonian $H$ of the form \eqref{Hquad} we can find $L \in \operatorname{Sp}(\R^{2n})$, such that $H \circ L$ satisfies assumptions of Lemma \ref{lem:h1}. In fact, by the H\"{o}rmander classification we know that positive definite matrices are symplectically diagonalizable. Therefore, we can assume that $H\circ L$ is of the form as in Lemma \ref{lem:h3}. Consequently, by Lemma \ref{lem:h1} and Lemma \ref{lem:h3} we get that $H \circ L$ is a tentacular Hamiltonian in the sense of Definition \ref{defTent}. By Remark \ref{rem:SpTent} we get that any Hamiltonian $H$ as in \eqref{Hquad} is a tentacular Hamiltonian and its zero level set is a tentacular hyperboloid.

\hfill $\square$

\section{The symplectically \texorpdfstring{$b$}{b}-compatible manifold}\label{sec_bcm}
In this section we will consider a symplectic manifold $T^*\R^n$ with the standard symplectic form $\omega_0$ and we will show how to construct for a $b$-symplectic manifold $(X,Z)$ symplectically $b$-compatible with $(T^*\R^n, \omega_0)$ so that the singular set $Z$ would correspond to the points\linebreak $\{q\in \R^n\ |\ \norm{q}=\infty\}$. We hope that this construction will give us information on how our original symplectic structure on $T^\ast \R^n$ behaves at infinity. This construction will work for an arbitrary $n \in \mathbb{N}$.

We will start by recalling the definition of \emph{adjunction spaces}. In general, an adjunction space is a topological space obtained by gluing two other topological spaces together via a homeomorphism. Further on we will introduce the $n$-spherical coordinates on $\R^n$ to replace the Cartesian coordinates. In other words, we will show that $\R^n\setminus \{0\}$ is diffeomorphic to $\R_+\times \s^{n-1}$ where $\R_+:=(0,\infty)$. This diffeomorphism lifts up to a symplectomorphism on the cotangent bundles. Using the spherical coordinates we define the McGehee transform. Finally, we construct the appropriate $b$-symplectic manifold compatible with $(T^*\R^n, \omega_0)$ as an adjunction space obtained by gluing two copies of $(T^*\R^n, \omega_0)$ together.
\subsection{Adjunction spaces}
The notion of \emph{adjunction spaces} was first introduced by Borsuk in \cite{borsuk1935} in 1935 and have been studied extensively since then. Generally, an adjunction space is a space obtained by gluing two topological spaces via a homeomorphism. However, in most cases in the literature (see for example \cite[Def. 6.1]{Dugundji1966}) the homeomorphism used in gluing the two topological spaces is defined on a closed set, whereas for our purposes we would prefer to glue them on an open set. The advantage of gluing two Hausdorff spaces via a homeomorphism on a closed set is that the resulted adjunction space is automatically Hausdorff. Unfortunately, the only two examples of an adjunction space defined via a homeomorphism on open set we have managed to find are recent papers by Luc and Placek \cite{luc2019} from 2019 and O'Connell \cite{oconnell2020} from 2020, where they consider non-Hausdorff spaces. Therefore, for the completeness of the argument, we decided to present the construction of adjunction spaces glued via a homeomorphism on an open set along with the proofs.

Hereby, we present the definition of an adjunction space we will use in this paper:
\begin{definition}\label{def_att}
Let $X$ and $Y$ be two topological spaces. Let $A\subseteq X$ be an open subset and let $f: A \to Y$ be an open embedding. We define
\[
X\cup_f Y:= X \sqcup Y /\sim,
\]
by identifying $a \sim f(a)$ for all $a\in A$. Moreover, we define topology on $X\cup_f Y$ by requiring the projection $\pi:X \sqcup Y \to X\cup_f Y$ to be continuous.
\end{definition}
Under the assumptions of the definition above we are able to prove the following properties of adjunctions spaces:
\begin{lemma}\label{lem_tpmf}
Let $X$ and $Y$ be two topological spaces. Let $A\subseteq X$ be an open subset and let $f: A \to Y$ be an open embedding. Then the projections $\pi|_X: X \to X\cup_f Y$ and $\pi|_Y:Y \to X\cup_f Y$ are open topological embeddings.
\end{lemma}
\begin{proof}
For an arbitrary, open subset $O\subset X$ we want to show that its image under the projection $\pi|_X$ is open in $X\cup_f Y$. Since $X\cup_f Y$ carries the quotient topology, this is equivalent to showing that its preimage under the projection $\pi$ is an open subset of $X\sqcup Y$. But by definition, we have:
\begin{equation*}
\pi^{-1}(\pi(O))=O\sqcup f(O\cap A).
\end{equation*}
Naturally, $O\cap A$ is open in $X$ and $f$ is an open embedding so $f(O\cap A)$ is open in $Y$. Thus $\pi^{-1}(\pi(O))$ is indeed an open subset of $X\sqcup Y$ and hence $\pi(O)$ is an open subset of $X\cup_f Y$. This proves that $\pi|_X$ is in fact an open map and thus a homeomorphism onto its image. The same is true for $\pi|_Y$ using a similar argument.
\end{proof}
\begin{lemma}\label{lem_smmf}
If $X$ and $Y$ are smooth manifolds, $f$ is a smooth embedding and the adjunction space $X\cup_f Y$ is Hausdorff, then there exists a unique smooth structure on $X\cup_f Y$, such that $X\cup_f Y$ is a smooth manifold and the projections $\pi|_X$ and $\pi|_Y$ are diffeomorphisms.
\end{lemma}
\begin{proof}
We will start with proving that $X\cup_f Y$ is a topological manifold. Since $X\cup_f Y$ is Hausdorff by assumption, it suffices for us to prove that it is also second countable and locally euclidean. By Lemma \ref{lem_tpmf} we know that both projections $\pi|_X$ and $\pi|_Y$ are open topological embeddings. Therefore\m{,} $\pi(X)$ and $\pi(Y)$ are locally euclidean as images of a homeomorphism of locally euclidean sets. Hence $X\cup_f Y$ is also locally euclidean since it is covered by locally euclidean open sets. The same proof also shows that $X\cup_f Y$ is second countable.

Now we will show how the smooth embeddings $\pi|_X$ and $\pi|_Y$ define a smooth structure on $X\cup_f Y$. 
More precisely, if we have a chart $(U,\phi)$ on $X$ then $\left(\pi(U),\phi\circ(\pi|_X)^{-1}\right)$ defines a chart on $X\cup_f Y$. Indeed, $\pi(U)$ being an open set and $\phi\circ(\pi|_X)^{-1}$ being a homeomorphism both follow from the fact that $\pi|_X$ is a homeomorphisms onto an open set in $X\cup_f Y$. Furthermore, one can easily see that this construction preserves the compatibility of charts.

We can use the same construction to transport the charts on $Y$ to charts on $X\cup_f Y$. Together these then cover all of $X\cup_f Y$. The charts coming from $X$ are in fact also compatible with the charts coming from $Y$ since $(\pi|_Y)^{-1}\circ\pi|_X=f$ is a diffeomorphism. Hence this defines a smooth atlas on $X\cup_f Y$.

The fact that $\pi|_X$ and $\pi|_Y$ are diffeomorphisms with respect to this smooth atlas follows from the construction. Lastly, we observe that any other smooth atlas with this property necessarily defines the same smooth structure.

\end{proof}

In a final step, we will show how one can use the same construction to glue together ($b$-)symplectic manifolds:
\begin{lemma}\label{lem_bglue}
Let $(M, \omega)$ be a symplectic manifold and let $(X, Z, \omega_b)$ be a $b$-symplectic manifold of the same dimension. Let $Y \subseteq X$ be an open subset disjoint with $Z$ and let $f: Y \to M$ be a symplectic embedding. If $M \cup_f X$ is a Hausdorff space then there exists a unique $b$-symplectic structure $\Omega_b$ on $(M \cup_f X, \pi(Z))$, such that $\pi^*(\Omega_b \big|_{\pi(X)})=\omega_b$ and $\pi^*(\Omega_b\big|_{\pi(M)})=\omega$.
\end{lemma}
\begin{proof}
Since symplectic embeddings are in particular smooth, $f$ is a diffeomorphism onto its image. By assumption the manifolds $M$ and $X$ have the same dimension, hence the image of $f$ is open in $M$ by the inverse function theorem. As such we can apply Lemma \ref{lem_tpmf} and \ref{lem_smmf} to get that $M \cup_f X$ is a smooth manifold. Furthermore, from the fact that $\pi|_X$ is a smooth (and open) embedding it follows that $\pi(Z)$ is a hypersurface of $M \cup_f X$. Hence we can look at $(M \cup_f X, \pi(Z))$ as a $b$-manifold.

We can use the equations $\pi^*(\Omega_b \big|_{\pi(X)})=\omega_b$ and $\pi^*(\Omega_b\big|_{\pi(M)})=\omega$ to define a $b$-2-form $\Omega_b$ on $(M \cup_f X, \pi(Z))$ . This is then well-defined since $f$ is a symplectomorphism onto its image by assumption. Lastly, it is non-degenerate as a $b$-2-form since $\omega$ and $\omega_b$ are.
\end{proof}
\subsection{The Spherical Coordinates}\label{sec_msc}

In this section we will define the $n$-spherical coordinate change on the cotangent bundle $T^\ast(\R^n\setminus\{0\})$.
As a starting point, we assume that we are working on $T^\ast \R^n$ with Cartesian coordinates $(q,p)$ and the standard symplectic form $\omega_0:=dq\wedge dp$. We will now look at $\s^{n-1}$ as the submanifold of $\R^n$ defined by
\[
\s^{n-1}=\{\psi\in\R^n\ \big|\ \psi^T\psi=1\},
\]
and denote $Q:= \R_+\times \s^{n-1}$. Then we can define the smooth map:
\begin{equation}\label{F}
\begin{array}{r c l}
F:\R^n\setminus\{0\} & \to & Q,\\
q & \mapsto & \left(\norm{q},\frac{q}{\norm{q}}\right).
\end{array}
\end{equation}
This is a smooth map with the smooth inverse:
\begin{equation}
\begin{array}{r c l}
F^{-1}:Q & \to & \R^n\setminus \{0\},\\
(r,\psi) & \mapsto & r\cdot\psi.
\end{array}
\end{equation}
Hence it is a diffeomorphism. We can now lift $F$ to get a symplectomorphism between the canonical symplectic forms on the cotangent bundles as in \cite[Cor.\ 2.2]{Silva2001}:
\begin{equation}\label{eq_lift}
\begin{array}{r c l}
F^\#:T^\ast(\R^n\setminus \{0\}) & \to & T^\ast Q,\\
(q,p) & \mapsto & \left(F(q),p\circ DF^{-1}\right).
\end{array}
\end{equation}
In the later part of this thesis we will do some explicit calculations, so we do need to fix some coordinates on $Q$. Fortunately, we can just take an arbitrary atlas on $\s^{n-1}$. By taking the product with the trivial chart on $\R_+$ we can modify this to get an atlas on $Q$. Now let $(O,\phi)$ be one of the charts from the arbitrary atlas on $\s^{n-1}$, i.e. $\phi:O\to U$ is a diffeomorphism, while $O\subset\s^{n-1}$ and $U\subset\R^{n-1}$ are some open subsets. We denote the inverse of $\phi$ by $\psi$. We will generally think of $\psi$ as an element of $Q$ dependent on $\phi$.

From now on we write $(r,P_r)$ for the coordinates on $T^*\R_+$ and $(\phi,\eta)$ for the coordinates on $T^\ast\s^{n-1}$. Without loss of generality we can usually assume that we are working inside $\R_+\times O$, the domain of the chart $Id\times\phi$. We can write the canonical symplectic form on $T^\ast Q$ using these new coordinates and by the symplectomorphism $F^\#$ it satisfies:
\begin{equation}
\left(F^\#\right)^*\left(dr\wedge dP_r+d\phi\wedge d\eta\right)=dq\wedge dp.
\end{equation}
The coordinate change for the base coordinates is given by the formula:
\[
q=F^{-1}\circ(Id\times\psi)(r,\phi)=r\cdot\psi(\phi).
\]
For the cotangent fiber coordinates we have the following lemma:
\begin{lemma}\label{lem_cc}
The coordinate change for the cotangent fiber coordinates is given by:
\begin{align}
p=&\,\left(P_r\frac{\partial r}{\partial q}^T+\frac{\partial \phi}{\partial q}^T\cdot\eta\right),\label{eq_cc1}
\\(P_r,\eta)=&\,\left(\frac{\partial q}{\partial r}^T\cdot p,\frac{\partial q}{\partial \phi}^T\cdot p\right).\label{eq_cc2}
\end{align}
\end{lemma}
\begin{proof}
In order to gain a better understanding of the map $F^\#$, we will first take a look at the map $DF^{-1}:TQ\to T(\R^n\setminus\{0\})$. This is naturally fiber-wise linear and as such can be represented by a matrix. In fact by \cite[Prop.\ 4.2.13]{Csikos} we have that with respect to the bases $\partial_q$ and $\left\lbrace\partial_r,\partial_\phi\right\rbrace$ the map $DF^{-1}$ is given by the matrix:
\[D((Id\times\phi)\circ F)^{-1}.\]
Now we notice that if we restrict $F^\#$ to a fiber, then this just corresponds to the dual map of $DF^{-1}$. By \cite[3.114]{axler1997} we know that the dual map with respect to the dual bases is just given by the transpose of the original matrix. Hence we get that with respect to the bases $\{dq\}$ and $\{dr,d\phi\}$ the map $F^\#$ restricted to a fiber is given by the matrix: \[\left(D((Id\times\phi)\circ F)^{-1}\right)^T.\]

In appropriate coordinates the map $(Id\times\phi)\circ F$ can be expressed as:
\[(Id\times\phi)\circ F(q)=(r,\phi)\,\in\,\R_+\!\times\R^{n-1}.\]
Using these coordinates we obtain the following expression for  the cotangent fiber coordinates:
\[
\begin{pmatrix}P_r\\\eta\end{pmatrix}=\left(D((Id\times\phi)\circ F)^{-1}\right)^T\cdot p=\begin{pmatrix}\frac{\partial q}{\partial r}&\frac{\partial q}{\partial \phi}\end{pmatrix}^T \cdot p=\begin{pmatrix}\frac{\partial q}{\partial r}^T\cdot p\\\frac{\partial q}{\partial \phi}^T\cdot p\end{pmatrix},
\]
which is exactly the equality \eqref{eq_cc2}. For \eqref{eq_cc1} we simply use the inverse matrix:
\begin{align*}
p=&\,\left(D((Id\times\phi)\circ F)\right)^T \cdot\begin{pmatrix}P_r\\\eta\end{pmatrix}=\begin{pmatrix}\frac{\partial r}{\partial q}\\\frac{\partial \phi}{\partial q}\end{pmatrix}^T\cdot\begin{pmatrix}P_r\\\eta\end{pmatrix}
\\=&\,\begin{pmatrix}\frac{\partial r}{\partial q}^T&\frac{\partial \phi}{\partial q}^T\end{pmatrix}\cdot\begin{pmatrix}P_r\\\eta\end{pmatrix}=P_r\frac{\partial r}{\partial q}^T+\frac{\partial \phi}{\partial q}^T\cdot\eta.
\end{align*}
\end{proof}
\begin{lemma}\label{lem_ps}
The two partial derivatives involving $r$ satisfy the following relation:
\begin{equation}\label{eq_pdqr}
\frac{\partial q_i}{\partial r}=\frac{\partial r}{\partial q_i}=\psi_i,
\end{equation}
where $i\in\{1,\dots, n\}$. As such they are both independent of $r$.
\end{lemma}
\begin{proof}
We can compute the first partial derivative using the formula $q=r\cdot\psi$ and the fact that $\psi$ only depends on $\phi$:
\[\frac{\partial  q_i}{\partial r}=\frac{\partial (r\cdot\psi_i)}{\partial r}=\psi_i.\]
For the second partial derivative we use the formula $r= \norm{q}=\sqrt{\sum_iq_i^2}$ instead:
 \[\frac{\partial r}{\partial q_i}=\frac{2q_i}{2\sqrt{\sum_iq_i^2}}=\frac{q_i}{r}=\psi_i.\]
\end{proof}
The partial derivatives involving $\phi$ are of course very difficult to compute, as we don't have an explicit formula for $\phi$. Instead, we will introduce the following notation:
\begin{equation}\label{def_mat}
\begin{array}{r c c l}
\text{We write} & \Psi & \text{for the matrix} & \Psi_{ij}:=r\frac{\partial \phi_j}{\partial q_i}\\
\text{and} &\Phi  &\text{for the matrix} & \Phi_{ij}:=\frac{1}{r}\frac{\partial q_j}{\partial \phi_i}.
\end{array}
\end{equation}
The following lemmas will allow us to derive the properties, which will be used in the calculations later on:
%The properties from the following lemmas will then be enough for the necessary calculations in this thesis.
\begin{lemma}\label{lem_in}
The matrices $\Phi$ and $\Psi$ are independent of $r$.
\end{lemma}
\begin{proof}
First we can prove that $\Phi$ is independent of $r$ by using the formula $q=r\cdot\psi$:
\[\Phi_{ij}=\frac{1}{r}\frac{\partial q_j}{\partial \phi_i}=\frac{1}{r}\frac{\partial (r\cdot\psi_{j})}{\partial \phi_i}=\frac{r}{r}\frac{\partial \psi_{j}}{\partial \phi_i}=\frac{\partial \psi_{j}}{\partial \phi_i}.\]
Since $\psi$ is only a function of $\phi$, this is independent of $r$. Next, we will proof:
\begin{equation}\label{eq_mi}
\begin{pmatrix}\psi^T\\\Psi^T\end{pmatrix}=\begin{pmatrix}\psi&\Phi^T\end{pmatrix}^{-1}.
\end{equation}
From this it will then follow that $\Psi$ is also independent of $r$. From what we have seen above we already know that the matrix $(\psi,\Phi^T)$ is completely independent of $r$. Hence we can conclude that its inverse, which includes $\Psi^T$ as a submatrix, is also independent of $r$.

So what is left to prove now is the equality \eqref{eq_mi}. First, observe that
that the Jacobian of inverse maps are inverses of each other:
 \[
\begin{pmatrix}\frac{\partial r}{\partial q}\hspace*{.1cm}\\\frac{\partial \phi}{\partial q}\end{pmatrix}=\,\begin{pmatrix}\frac{\partial q}{\partial r}&\frac{\partial q}{\partial \phi}\end{pmatrix}^{-1}.
\]
Now, using the equation above, the notation from \eqref{def_mat} and the result from Lemma \ref{lem_ps} we can calculate:
\begin{align*}
\begin{pmatrix}\psi^T\\\Psi^T\end{pmatrix} & =\begin{pmatrix}\hspace*{.1cm}\frac{\partial r}{\partial q}\hspace*{.1cm}\\r\cdot\frac{\partial \phi}{\partial q}\end{pmatrix}=\begin{pmatrix}1&0\\0&r\cdot Id\end{pmatrix}
\cdot\begin{pmatrix}\hspace*{.1cm}\frac{\partial r}{\partial q}\hspace*{.1cm}\\\frac{\partial \phi}{\partial q}\end{pmatrix}\\
&=\left(\begin{pmatrix}\frac{\partial q}{\partial r}&\frac{\partial q}{\partial \phi}\end{pmatrix}
\cdot\begin{pmatrix}1&0\\0&\frac{1}{r}\cdot Id\end{pmatrix}\right)^{-1}\\
& =\begin{pmatrix}\frac{\partial  q}{\partial r}&\frac{1}{r}\cdot\frac{\partial q}{\partial \phi}\end{pmatrix}^{-1}=\begin{pmatrix}\psi&\Phi^T\end{pmatrix}^{-1}.
\end{align*}
This proves equation \eqref{eq_mi}. Since both $\psi$ and $\Phi$ are independent of $r$, then by \eqref{eq_mi} so is $\Psi$.
\end{proof}
\begin{lemma}\label{lem_for}
Additionally, the following equalities are satisfied:
\begin{align}
\Psi^T\psi=&\,0,\label{eq_for1}
\\\Phi\psi=&\,0,\label{eq_for2}
\\\Phi\Psi=&\,Id,\label{eq_for3}
\\\Psi\Phi+\psi\psi^T=&\,Id.\label{eq_for4}
\end{align}
\end{lemma}
\begin{proof}
From \eqref{eq_mi} we have:
\[Id=\begin{pmatrix}\psi^T\\\Psi^T\end{pmatrix}\begin{pmatrix}\psi&\Phi^T\end{pmatrix}=\begin{pmatrix}\psi^T\psi&\psi^T\Phi^T\\\Psi^T\psi&\Psi^T\Phi^T\end{pmatrix}.\]
This immediately gives us the equalities \eqref{eq_for1}, \eqref{eq_for2} and \eqref{eq_for3}. For the last equality \eqref{eq_for4}, we can apply \eqref{eq_mi} the other way around to get:
\[Id=\begin{pmatrix}\psi&\Phi^T\end{pmatrix}\begin{pmatrix}\psi^T\\\Psi^T\end{pmatrix}=\psi\psi^T+\Phi^T\Psi^T,\]
which is of course exactly the transpose of \eqref{eq_for4}.
\end{proof}
Substituting the definitions of $\Psi$ and $\Phi$ from Definition \ref{def_mat} into \eqref{eq_cc1} and \eqref{eq_cc2} gives us:
\begin{align}
 p=&\,P_r\psi+\frac{1}{r}\Psi\cdot\eta\label{eq_cc3},
\\(P_r,\eta)=&\,(\psi^T  p ,r\Phi p)\label{eq_cc4}.
\end{align}
Where we know that $\psi$, $\Phi$ and $\Psi$ are all independent of $r$ by Lemma \ref{lem_in}.

\subsection{The McGehee transform}\label{sec_McG}

The aim of this subsection is to construct a $b$-symplectic manifold symplectically $b$-compatible with $(T^*(\R^n \setminus \{0\}), \omega_0)$. The construction itself is not new and has been used previously to study the restricted $3$-body problem for $n=2$ in \cite{Delshams2019, Kiesenhofer2017, miranda2021}. Below we present its generalization to an arbitrary dimension $n>1$.

To keep things concise we introduce the following notation:
\begin{equation}\label{eq_defWY}
\begin{array}{r c l}
\widetilde{W} & := & \R^2\times T^*\s^{n-1,}\\
\widetilde{Y} & := & \{0\}\times\R\times T^*\s^{n-1},\\
W^+ & := & (0, +\infty)\times\R\times T^*\s^{n-1},\\
W^- & := & (-\infty,0)\times\R\times T^*\s^{n-1}. 
\end{array}
\end{equation}
This way $\widetilde{Y}$ is a hypersurface in $\widetilde{W}$, which makes $(\widetilde{W},\widetilde{Y})$ into a $b$-manifold. Moreover, we have the decomposition $\widetilde{W}= W^-\sqcup \widetilde{Y}\sqcup  W^+$. Note that $W^+$ and $W^-$ are both diffeomorphic to $T^*Q$. Hence the $b$-manifold $(\widetilde{W},\widetilde{Y})$ is $b$-compatible with $T^*Q$. In fact, by the diffeomorphism $F^\#:T^\ast(\R^n\setminus\{0\})\to T^\ast Q$ from \eqref{eq_lift} the $b$-manifold $(\widetilde{W},\widetilde{Y})$ is $b$-compatible with $T^*(\R^n \setminus \{0\})$.

Our aim now will be to construct a $b$-symplectic structure on $(\widetilde{W},\widetilde{Y})$, such that it would be symplectically $b$-compatible with $(T^*(\R^n \setminus \{0\}), \omega_0)$. Later on, in the next subsection, we will glue in the missing points in order to construct a $b$-symplectic manifold $(W,Y, \omega_b)$ which would be symplectically $b$-compatible with the whole $(T^*\R^n, \omega_0)$.

Next we want to apply the McGehee transformation which is often used to study the restricted three body problem \cite{miranda2021,}. It was first introduced by McGehee to study parabolic orbits \cite{mcgehee1973}. We can generalize this transformation to our setting by defining the following map:
\begin{equation}\label{eq_tau}
\begin{array}{r c l}
\tau\mcg: \widetilde{W}=\R^2 \times T^*\s^{n-1} & \to & T^*\R_+ \times T^* \s^{n-1}=Q,\\
(z, P_r, \phi, \eta) & \mapsto & \left( \frac{2}{z^2}, P_r, \phi, \eta \right).
\end{array}
\end{equation}
We will denote its restrictions to $W^-$ and $W^+$ by $\tau\mcg^-$ and $\tau\mcg^+$ respectively. Observe that both $\tau\mcg^-$ and $\tau\mcg^+$ are diffeomorphisms between $W^-, W^+$ and $T^*Q$, respectively. Now we will use these diffeomorphisms to pull back the standard symplectic form $\omega_0$ from $T^*Q$ to $\widetilde{W}\setminus \widetilde{Y}$.

We can calculate this pull back to be:
\begin{equation}\label{eq1}
\left(\tau\mcg^\pm\right)^*(\omega)=-\frac{4}{z^3}dz\wedge dP_r+d\phi\wedge d\eta.
\end{equation}
 Note that the symplectic form from \eqref{eq1} extends smoothly to a $b^3$-2-form on $\widetilde{W}$. Hence we make the obvious definition:
\begin{equation}\label{eq2}
\omega\mcg:=-\frac{4}{z^3}dz\wedge dP_r+d\phi\wedge d\eta.
\end{equation}
\begin{lemma}
The $b$-manifold $(\widetilde{W},\widetilde{Y})$ together with the $b^3$-form $\omega\mcg$ defines a $b^3$-symplectic manifold. Furthermore, $(\widetilde{W},\widetilde{Y},\omega\mcg)$ is symplectically $b$-compatible with $T^\ast(\R^n\setminus\{0\})$ with the standard symplectic form.
\end{lemma}
\begin{proof}
The $b^3$-form $\omega\mcg$ is clearly closed and non-degenerate as a $b^3$-form, hence $(\widetilde{W},\widetilde{Y},\omega\mcg)$ is in fact a $b^3$-symplectic manifold.

Since we used pull backs to define $\omega\mcg$, $\tau_{McG}^+$ is clearly a symplectomorphism between $W+$ and $T^*Q$. The same is true for $\tau_{McG}^-$ and $W^-$ respectively. So we have that $\widetilde{W}$ is symplectically $b$-compatible with $T^*Q$. But we have seen in Section \ref{sec_msc} that $F^\#$ is a symplectomorphism between $T^\ast( \R^n\setminus\{0\})$ and $T^\ast Q$. Hence $\widetilde{W}$ is also symplectically\linebreak $b$-compatible with $T^\ast( \R^n\setminus\{0\})$ with the standard symplectic form.
\end{proof}

\subsection{The symplectically \texorpdfstring{$b$}{b}-compatible manifold}\label{sec_glm}

The aim of this subsection will be to modify $(\widetilde{W},\widetilde{Y},\omega\mcg)$ in a certain way as to make it symplectically $b$-compatible with all of $(T^*\R^n, \omega_0)$. To do that we define
\[W:=((\R^{2n}\times\{- 1,1\})\sqcup\widetilde{W})/_\sim,\]
where the relation $\sim$ is defined by:
\begin{equation}\label{eq_rel}
\forall\ (x,y) \in W^\pm \times \R^{2n} \qquad x \sim (y, \pm 1) \quad \Leftrightarrow \quad \tau\mcg^\pm(x)=F^\#(y),
\end{equation} 
where $F^\#$ is the symplectomorphism as in \eqref{eq_lift}. If we write $\pi$ for the projection\linebreak $\pi:(\R^{2n}\times\{- 1,1\})\sqcup\widetilde{W}\to W$, then we can define the topology of $W$ by requiring $\pi$ to be continuous.

Equivalently, we can use Definition \ref{def_att} to define $W$ as
\begin{equation}\label{eq_defW}
W:= (\R^{2n}\times\{- 1,1\})\cup_f \widetilde{W},\qquad \text{where}\qquad f: \widetilde{W}\setminus \widetilde{Y} \to T^*\R^n\times\{-1,1\},
\end{equation}
is an open embedding %for all $x\in T^*(\R^n\setminus \{0\})$ 
defined by
\begin{equation}\label{eq_deff}
f(x):= \begin{cases}
\big(\left(F^\#\right)^{-1}\circ \tau\mcg^-(x),-1\big) & \text{for}\ x\in W^-,\\
\big(\left(F^\#\right)^{-1}\circ \tau\mcg^+(x),1\big) & \text{for}\ x\in W^+.
\end{cases}
\end{equation}
 The image of $f$ is $T^*(\R^n\setminus \{0\})\times \{-1,1\}$ and $f$ is a symplectic embedding, i.e. $f^*\omega_0= \omega\mcg$.
 
%Our aim now is to use Lemma \ref{lem_bglue} to show that $W$ is a smooth manifold. To to that, we have to show that the topology on $W$ induced by $\pi$ is Hausdorff.
Denote $Y:=\pi(\widetilde{Y})$. Our aim now is to use Lemma \ref{lem_bglue} to show that $(W,Y)$ is a smooth $b$-manifold that can be equipped with a unique $b$-symplectic structure $\omega_b$, such that 
$$
\pi^*(\omega_b\big|_{\pi(\R^{2n}\times\{- 1,1\})})=\omega_0
\qquad \text{and} \qquad
\pi^*(\omega_b\big|_{\pi(\widetilde{W})})=\omega\mcg.
$$
To use Lemma \ref{lem_bglue} we first need to show that the topology on $W$ induced by $\pi$ is Hausdorff.

\begin{lemma}\label{lem_xschd}
The topological space {$W$} is Hausdorff.
\end{lemma}
\begin{proof}
By Lemma \ref{lem_tpmf} we know that the projection $\pi:(\R^{2n}\times\{- 1,1\})\sqcup\widetilde{W}\to W$ is an open map. Due to the fact that $(\R^{2n}\times \{-1, 1\})\sqcup\widetilde{W}$ is Hausdorff, it follows that the two sets $\pi(\R^{2n}\times\{- 1,1\})$ and $\pi(\widetilde{W})$ are open and Hausdorff subsets of $W$. Since they cover $W$, it only remains to be shown that their complements can be separated with open neighborhoods.

Observe that 
\[
\pi(\widetilde{W}) =W \setminus\pi( T^*_0\R^n\times\!\{-1,1\}),\qquad\text{and} \qquad
\pi(\R^{2n}\times\{- 1,1\}) =W\setminus Y.
\]
 Now we define the open sets 
\[
\mathcal{O}:= \left\lbrace (q,p)\in T^*\R^n\ \Big|\ \norm{q}<1 \right\rbrace \times \{-1,1\}, \qquad
\mathcal{V}:=(-1,1)\times\R \times T^*\s^{n-1}\subseteq  \widetilde{W}.
\]
Since $\pi$ is an open map so $\pi(O)$ and $\pi(V)$ are open subsets of $\widetilde{W}$. Additionally, we have the following inclusions
\begin{align*}
W\setminus \pi(\widetilde{W})= \pi(T^*_0\R^n\times\!\{-1,1\}) \subseteq \pi (\mathcal{O})  & \subseteq \pi(\R^{2n}\times\{- 1,1\}) =W\setminus Y,\\
W \setminus \pi(\R^{2n}\times\{- 1,1\})= Y \subseteq \pi(\mathcal{V}) & \subseteq \pi(\widetilde{W}) =W\setminus\pi(T^*_0\R^n\times\!\{-1,1\}).
\end{align*}

On the other hand, by \eqref{eq_lift}, \eqref{eq_tau} and \eqref{eq_deff} we have
\begin{align*}
f(\mathcal{V}\setminus \widetilde{Y}) & = \left( F^\#\right)^{-1}\circ \tau\mcg \left( (0,1)\times\R \times T^*\s^{n-1}\right) \times\{-1,1\}\\
& = \left(F^\#\right)^{-1} \left( (2, +\infty)\times \R  \times T^*\s^{n-1} \right)\times\{-1,1\}\\
& = \left\lbrace (q,p)\in T^*\R^n\ \big|\ \norm{q}>2 \right\rbrace \times\{-1,1\}
\end{align*}
Consequently, the sets $f(\mathcal{V}\setminus \widetilde{Y})$ and $\mathcal{O}$ are disjoint in $T^*\R^n\times\{-1,1\}$ and therefore by \eqref{eq_defW} the sets $\pi(\mathcal{O})$ and $\pi(\mathcal{V})$ are disjoint in $W$.

We can now finally prove that $W$ is Hausdorff. Consider two distinct points $x,y\in X$. If both are contained in $\pi(\widetilde{W})$ then we are done, since we know $\pi(\widetilde{W})$ to be Hausdorff. The same is true if both are contained in $\pi(\R^{2n}\times\{- 1,1\})$. The last remaining possibility is that $x\in W\setminus\pi( \widetilde{W})$ and $y\in W\setminus\pi(\R^{2n}\times\{-1, 1\})$ (or the other way around). In this case they are separated by the open and disjoint neighbourhoods $\pi(\mathcal{O})$ and $\pi(\mathcal{V})$.
\end{proof}
\begin{proposition}\label{prop_bsm}
The pair $(W,Y)$ is a $b^3$-manifold and it has a $b^3$-2-form $\omega_b$ such that $(W, Y,\omega_b)$ is a $b^3$-symplectic manifold and is symplectically $b$-compatible with $T^*\R^n$ carrying the standard symplectic structure.
\end{proposition}
\begin{proof}
By Lemma \ref{lem_xschd} we have that the topology on $W= (\R^{2n}\times\{- 1,1\})\cup_f \widetilde{W}$ is Hausdorf. Consequently, $W$ satisfies all the assumptions of Lemma \ref{lem_bglue}. Hence we immediately get that $( W, Y)$ is in fact a $b$-manifold. Furthermore, since the function $f$ is a symplectic embedding of $\widetilde{W}\setminus \widetilde{Y}$ into $T^*\R^n \times\{-1,1\}$, we can infer by Lemma \ref{lem_xschd} that $( W, Y)$
carries a unique $b^3$-symplectic form $\omega_b$ such that the maps $\pi|_{T^*\R^n\times \{- 1,1\}}$ and $\pi|_{\widetilde{W}}$ are\linebreak $b^3$-symplectomorphisms onto their images.

Observe that $\pi(T^*\R^n\times\{-1,1\})=W\setminus Y$. In fact, by construction we have that\linebreak $T^*\R^n\times\{- 1,1\}$ with the standard symplectic structure is symplectomorphic to $W\setminus Y$. Naturally the trivial double-cover $T^*\R^n\times\{- 1, 1\}\to T^*\R^n$ is a symplectomorphism when restricted to either connected component. So it follows that each connected component of $W\setminus Y$ is symplectomorphic to $(T^*\R^n, \omega_0)$, which is what we wanted to prove.
\end{proof}
\begin{remark}
Note that this in general might not be the only $b^3$-symplectic structure on $(W,Y)$ such that it is $b$-symplectically compatible with $T^*\R^n$. In Appendix \ref{app_nonun} we show an example of two non-$b$-symplectomorphic $b^m$-symplectic structures which become symplectomorphic once you remove the singular hypersurface. As a consequence the corresponding $b^m$-symplectic manifolds are symplectically $b$-compatible with the same symplectic manifolds.
\end{remark}
\subsection{Liouville \texorpdfstring{$b$}{b}-vector field}
Our goal is to prove Theorem \ref{thm_main}. One way of proving the existence of a contact structure on a given hypersurface is by finding a Liouville vector field, which is everywhere transverse to this hypersurface. In this subsection we will construct a Liouville $b^3$-vector field on $(W,Y, \omega_b)$, which we will use later on in the next section to establish the contact structure.

Let $h: \R_+\to (0,2]$ be a smooth, positive function, such that
\[
h(t):= \begin{cases} 1 & \text{for}\ t \leq \frac{1}{2},\\
\frac{1}{t} & \text{for}\ t\geq 1.
\end{cases}
\]
Define a vector field $V$ on $T^*\R^n$ in the following way:
\begin{equation}\label{eq_defVhat}
\widehat{V}(q,p):=(p+h(\norm{q})q)\partial_p\qquad (q,p)\in T^*\R^n.
\end{equation}
Note that $V$ is smooth since $h$ is constant on a neighborhood of the set $\{\norm{q}=0\}$. Furthermore, $\widehat{V}$ satisfies $d\iota_{\widehat{V}} \omega_0=\omega_0$. In other words, $\widehat{V}$ is a Liouville vector field on $(T^*\R^n, \omega_0)$.

\begin{remark}
By $p\partial_p$ we of course mean $\sum_{i=1}^np_i\partial_{p_i}$. In order to interpret this as a matrix multiplication it should hence be written as $p^T\partial_p$. For the sake of simplicity sake we will stick to the former choice of notation whenever reasonable.
\end{remark}

\begin{lemma}
The Liouville vector field $\widehat{V}$ on $(T^*\R^n, \omega_0)$ extends uniquely to a  Liouville $b^3$-vector field on $(W, Y, \omega_b)$.
\end{lemma}

\begin{proof}
The function $f$ from \eqref{eq_deff} is a diffeomorphism onto its image. With a slight abuse of notation we will denote its inverse from its image by $f^{-1}: T^*(\R^n\setminus \{0\})\times\{-1,1\}\to \widetilde{W}\setminus \widetilde{Y}$.
Since $\widehat{V}$ is a Liouville vector field on $(T^*\R^n, \omega_0)$, to show that it extends uniquely to a  Liouville $b^3$-vector field on $(W, Y, \omega_b)$ it suffices to show that $df^{-1}(\widehat{V})$ extends to a $b^3$-vector field on $(\widetilde{W}, \widetilde{Y}, \omega\mcg)$. To prove that we will start with calculating $dF^\#(\widehat{V})$.

Note that $r$ and $\psi$ are functions of $q$ and are independent of $p$. Therefore, by \eqref{eq_cc4} we have
\[
\partial_p=\frac{\partial P_r}{\partial p}\partial_{P_r}+ \frac{\partial \eta}{\partial p}\partial_\eta= \psi^T\partial_{P_r}+r \Phi \partial_\eta.
\]
Plugging that into \eqref{eq_defV} and using \eqref{eq_for1}, \eqref{eq_for2}, \eqref{eq_for3}, \eqref{eq_cc3}, \eqref{eq_cc4} and the fact that $\psi^T \psi=1$ we obtain
\[
dF^\# (\widehat{V}) = \left( P_r \psi^T + \frac{1}{r}\eta^T\Psi^T + h(r) r \psi^T \right)\left(\psi\partial_{P_r}+r \Phi^T\partial_\eta\right)=\left(P_r +r h(r)\right)\partial_{P_r}+\eta \partial_\eta.
\]
Consequently,
\[
df^{-1}(\widehat{V})=\left(P_r +\frac{2}{z^2}h\left(\frac{2}{z^2}\right)\right)\partial_{P_r}+\eta \partial_\eta,
\]
where $\frac{2}{z^2}h(\frac{2}{z^2})$ is a smooth, positive function on $\R$, such that
\[
\frac{2}{z^2} h\left(\frac{2}{z^2}\right)= \begin{cases}
1 & \text{for}\ |z|\leq \sqrt{2},\\
\frac{2}{z^2} & \text{for}\ |z|\geq 2. 
\end{cases}
\]
Since $d f^{-1}(\widehat{V})$ does not depend on $\partial_z$ or $z$ in a neighbourhood of $\widetilde{Y}=\{0\}\times\R \times T^*\s^{n-1}$ it follows that it can be uniquely extended to a $b^3$-vector field $\widetilde{V}$ on $(\widetilde{W},\widetilde{Y})$. A straightforward calculation shows that $d\iota_{\widetilde{V}}\omega\mcg=\omega\mcg$, which proves that $\widetilde{V}$ is a Liouville $b^3$-vector field on $(\widetilde{W},\widetilde{Y}, \omega\mcg)$.

We define a vector field $V$ on $(W,Y)$ by setting:
\begin{equation}\label{eq_defV}
V(w):=\begin{cases}
d\pi(\widehat{V}) & \text{for}\ w \in \pi(T^*\R^n\times\{-1,1\}),\\
d\pi(\widetilde{V}) & \text{for}\ w \in \pi(\widetilde{W}).
\end{cases}
\end{equation}
Since $\widehat{V}$ is a Liouville vector field on $T^*\R^n$ and $\widetilde{V}$ is a Liouville $b^3$-vector field on $(\widetilde{W},\widetilde{Y},\omega\mcg)$ satisfying $ df^{-1}(\widehat{V})= \widetilde{V}$, therefore by \eqref{eq_deff} $V$ is a well defined $b^3$-vector field on $(W, Y)$. On the other hand, by construction $\pi: (T^*\R^n\times\{-1,1\}) \sqcup \widetilde{W} \to W$ restricted to $T^*\R^n\times\{-1,1\}$ and $\widetilde{W}$ separately are symplectomorphisms. Consequently, $V$ is a Liouville $b^3$-vector field on $(W, Y, \omega_b)$.
\end{proof}
\section{The \texorpdfstring{$b$}{b}-contact structure}\label{sec:bcont}
The aim of this section will be to prove Theorem \ref{thm_main}. However, rather than using a Hamiltonian of the form as in \eqref{defH} we will consider a Hamiltonian $H: M \times T^*\R^m \to \R$ defined
\begin{equation}\label{eq_H}
H(x, q, p):= H_0(x)+p^T B q, \qquad x\in M,\ (q,p)\in T^*\R^k,
\end{equation}
where the function $H_0:M \to \R$ satisfies \eqref{defH0}, however $B$ is a non-degenerate, $k \times k$ matrix, such that $B+B^T$ is positive definite. By Proposition \ref{prop_JAB} the Hamiltonians of the form as in \eqref{defH} and in \eqref{eq_H} are equivalent up to a symplectic change of coordinates. In particular, their zero level sets define the same symplectic hyperboloid. 
Therefore, without loss of generality, from now on we will consider Hamiltonians of the form as in \eqref{eq_H}.
We denote by $S\subseteq M \times T^*\R^k$ the zero level set of $H$. %, i.e. our symplectic hyperboloid.
On the other hand we define a $b$-symplectic manifold $(X,Z, \omega \oplus \omega_b)$, where
\begin{equation}\label{eq_X}
X:= M \times W,\qquad \text{and}\qquad Z:= M \times Y,
\end{equation}
and $(W,Y, \omega_b)$ is the $b$-symplectic manifold defined in previous section symplectically\linebreak $b$-compatible with $(T^*\R^k, \omega_0)$. Denote $\Pi:=id \oplus \pi$, where $\pi$ is the projection\linebreak $\pi: T^*\R^k \times\{-1,1\}\sqcup \widetilde{W}\to W$. By Remark \ref{rem_prod} and Proposition \ref{prop_bsm} the $b$-manifold $(X,Z, \omega \oplus \omega_b)$ is symplectically $b$-compatible with $(M \times T^*\R^k, \omega \oplus \omega_0)$.

Note that we have the following inclusions
\begin{align}
S\times\{-1,1\} & \subseteq M\times T^*\R^k\times\{-1,1\},\nonumber \\
\Pi(S \times\{-1,1\}) & \subseteq X\setminus Z.\label{eq_St}
\end{align}
By Definition \ref{def_img} the $b$-image of $S$ in $( X, Z, \omega \oplus \omega_b)$ is equal to
\begin{equation}\label{eq_sb}
S_b:=\overline{\Pi (S\times\{-1,1\})}.
\end{equation}
Our first goal will be to show that $S_b$ is a smooth hypersurface in $X$.

\begin{lemma}\label{lem_Sb}
Let $S\subseteq M\times T^*\R^k$ be a hypersurface defined as a zero level set of a Hamiltonian defined in \eqref{eq_H}. Let $(X,Z, \omega \oplus \omega_b)$ be a $b$-symplectic manifold defined in \eqref{eq_X}.
Then the $b$-image of $S$ in $X$ defined as in \eqref{eq_sb} is a smooth hypersurface.
\end{lemma}

\begin{proof}
Denote $\widetilde{X}:= M \times\widetilde{W}$. Note that $\{X\setminus Z, \Pi(\widetilde{X})\}$ is an open cover of $X$. To show that $S_b$ is a smooth hypersurface in $X$, if suffices to show that its restrictions to both sets $X\setminus Z$ and $\Pi(\widetilde{X})$ are smooth. 
By construction, the restriction of the projection $\Pi$ to the set $ M \times T^*\R^k \times\{-1,1\}$ is a diffeomorphism onto its image, hence the set $ \Pi(S\times\{-1,1\})$ is a smooth hypersurface in $X\setminus Z$. Following from \eqref{eq_sb} and the fact that $S$ is a closed set we have $S_b\cap (X\setminus Z)= \Pi(S\times\{-1,1\})$ and hence it is also smooth. Therefore, what is left to be shown is that the restriction of $S_b$ to $\Pi(\widetilde{X})$ is a smooth hypersurface.

We will start with expressing $S_b$ in the coordinate chart of $\Pi(\widetilde{X})$. By \eqref{eq_St} we know that $S_b\cap(X\setminus Z)=\Pi (S\times\{-1,1\})$, so that its boundary is a subset of the critical set $Z$. On the other hand, by Lemma \ref{lem_xschd} we know that $\pi(T^*_0\R^k\times\{-1,1\})$ and $Y$ are separated by open sets in $W$. Therefore to determine $S_b\cap Z$ it suffices to analyze $S_b$ in the coordinates around $Z$. In other words, we will restrict ourselves to $S_b \setminus \Pi\left(M \times T^*_0\R^k\times\{-1,1\}\right)=S_b\cap \Pi(\widetilde{X})$. 

By \eqref{eq_rel} we obtain
\begin{align*}
\Pi (S\times\{-1,1\}) \cap \Pi(\widetilde{X})  & =\Pi ((S\cap (M \times T^*(\R^k\setminus\{0\}))\times\{-1,1\}))\\
%& = \Pi (\tilde{S}\times\{-1,1\})\\
%& = \Pi \left(\tau\mcg)^{-1}(F^\#(\tilde{S}))\right)
& = \Pi \left((id \oplus \tau\mcg)^{-1}(id \oplus F^\# (S\cap (M\times T^*(\R^k\setminus\{0\})))\right)\\
& =  \Pi \left((id \oplus \tau\mcg)^{-1}(id \oplus F^\# (H^{-1}(0)\cap (M\times T^*(\R^k\setminus\{0\})))\right)\\
& = \Pi \left(\left(H\circ (id \oplus (F^\#)^{-1}\circ\tau\mcg)\right)^{-1}(0)\right)\\
& = \Pi (S\mcg),
\end{align*}
where we denote $S\mcg:=\left(H\circ (id \oplus (F^\#)^{-1}\circ\tau\mcg)\right)^{-1}(0)\subseteq M\times Q$.

By construction the restriction of the projection $\Pi|_{M\times T^*\R^k\times\{-1,1\}}$ and $\Pi|_{\widetilde{X}}$ are  homeomorphisms onto their image, therefore we have
\[
S_b=\Pi \left(S\times\{-1,1\}\sqcup \overline{S\mcg}\right).
\]
To analyze $S\mcg$ we first express $H$ in spherical coordinates $(r,P_r,\phi,\eta)\in T^*\R_+ \times T^*\s^{k-1}$ using \eqref{eq_cc1}:
\begin{align*}
H \circ \left(id \oplus (F^\#)^{-1}\right)(x,r,P_r,\phi,\eta)&= H_0(x)+(r\psi)^TB\left(P_r\psi+\frac{1}{r}\Psi\eta\right)
\\&= H_0(x)+rP_r\psi^TB\psi+\psi^TB\Psi\eta,
\end{align*}
where we keep in mind that by Lemma \ref{lem_in} $\psi$ and $\Psi$ are functions of $\phi \in \s^{ k-1}$. Consequently the zero level set is given by the equation:
\[
\left(H \circ \left(id \oplus (F^\#)^{-1}\right)\right)^{-1}(0)=\left\{H_0(x)+rP_r\psi^TB\psi+\psi^TB\Psi\eta=0\right\}\subseteq M\times Q.
\]
Now we can apply the McGehee Transformation and look at the preimage in $\widetilde{X}\setminus \widetilde{Z}$, where $\widetilde{Z}:= M\times Y$.
By substituting $r=\frac{2}{z^2}$ into the equation above, we obtain
\begin{align}
S\mcg &= \left\{H_0(x)+\frac{2}{z^2}P_r\psi^TB\psi+\psi^TB\Psi\eta=0\right\}\nonumber\\
 & =\left\{2P_r\psi^TB\psi+z^2\left(H_0(x)+\psi^TB\Psi\eta\right)=0\right\}{\subseteq \widetilde{X}\setminus \widetilde{Z}}\label{eq_hyp2}
\end{align}
In order to write this more succinctly, we define the following smooth map on $\widetilde{X}$:
\begin{equation}\label{eq_Hmcg}
H\mcg(x,z,P_r,\phi,\eta):=2P_r\psi^TB\psi+z^2\left(H_0(x)+\psi^TB\Psi\eta\right).
\end{equation}
With this we can reformulate \eqref{eq_hyp2} to get:
\begin{equation}\label{eq_hyp3}
S\mcg=H\mcg^{-1}(0)\cap (\widetilde{X}\setminus \widetilde{Z}).
\end{equation}
Now we claim that the closure of $S\mcg$ in $\widetilde{X}$ is exactly $H\mcg^{-1}(0)$ which would then give us:
\begin{equation}\label{eq_cl}
\overline{S\mcg}=H\mcg^{-1}(0).
\end{equation}

First we will prove:
\begin{equation}\label{eq_hyp4}
\overline{S\mcg}\subseteq H\mcg^{-1}(0).
\end{equation}
Since $H\mcg^{-1}(0)$ is closed in $\widetilde{X}$ and it contains $S\mcg$ by \eqref{eq_hyp3}, it immediately follows that it also contains its closure, which proves \eqref{eq_hyp4}.

Next, we want to show:
\begin{equation}\label{eq_hyp5}
H\mcg^{-1}(0)\subseteq \overline{S\mcg}.
\end{equation}
By definition $S\mcg$ is a zero level set, so it is a closed set in $\widetilde{X}\setminus \widetilde{Z}$. Because taking the closure in $\widetilde{X}\setminus \widetilde{Z}$ is the same thing as restricting the closure in $\widetilde{X}$ to  $\widetilde{X}\setminus \widetilde{Z}$, it follows that $S\mcg=\overline{S\mcg}\cap(\widetilde{X}\setminus \widetilde{Z})$. Together with \eqref{eq_hyp3} we then have:
\begin{equation}\label{eq_hyp6}
S\mcg\cap(\widetilde{X}\setminus \widetilde{Z})=H\mcg^{-1}(0)\cap(\widetilde{X}\setminus \widetilde{Z}).
\end{equation}
Thus it only remains to be shown that:
\begin{equation}\label{eq_hyp7}
H\mcg^{-1}(0)\cap \widetilde{Z}\subseteq \overline{S\mcg}.
\end{equation}
Equivalently, we need to show  that we can reach all the points in $H\mcg^{-1}(0)\cap \widetilde{Z}$ with sequences in $S\mcg$. To figure out what $H\mcg^{-1}(0)\cap \widetilde{Z}$ is, we substitute $z=0$ in the equation $H\mcg(x, z,P_r,\phi,\eta)=0$. This gives us:
\[H\mcg^{-1}(0)\cap \widetilde{Z}=\{2P_r\psi^TB\psi=0,z=0\}.\]
We can rewrite $2\psi^TB\psi=\psi^T\left(B+B^T\right)\psi$. Since we know that $B+B^T$ is positive definite and $\psi$ is a length 1 vector, this is always positive. %So the set above can be simplified to $\{P_r=0,z=0\}$.
Consequently, we have
\[H\mcg^{-1}(0)\cap \widetilde{Z}=\{P_r=0,z=0\}.\]

Now let's fix a point $w_0:=(x_0,z_0,P_{0},\psi_0,\eta_0) \in H\mcg^{-1}(0)\cap  \widetilde{Z}$. Of course, this immediately gives us $z_0=P_{0}=0$. We can define a sequence of points $w_n:=(x_n,z_n,P_{n},\psi_n,\eta_n)\in \widetilde{W}$ by setting $x_n:=x_0,z_n:=\frac{1}{n},\psi_n:=\psi_0$ and $\eta_n:=\eta_0$. We can then choose $P_n$ such that $w_n$ lies in $S\mcg$. By \eqref{eq_hyp3} we know that this is equivalent to:
\begin{align*}
0=H\mcg(x_n,z_n,P_n,\psi_n,\eta_n) & = 2P_n\psi_0^TB\psi_0+\frac{1}{n^2}\left(H_0(x_0)+\psi_0^TB\Psi\eta_0\right),\\
P_n & = - \frac{1}{n^2}\cdot\frac{H_0(x_0)+\psi_0^TB\Psi\eta_0}{2\psi_0^TB\psi_0}.
\end{align*}
Since $2\psi_0^TB\psi_0$ is positive as above, this is just a constant times $\frac{1}{n^2}$. So both $z_n$ and $P_{n}$ approach 0 (and thus $z_0$ and $P_{0}$ respectively) for $n\to\infty$ and we have:
\[\lim_{n\to\infty}(x_n,z_n,P_{n},\psi_n,\eta_n)=(x_0,z_0,P_{0},\psi_0,\eta_0).\]
This means that $w_0$ is indeed part of $\overline{S\mcg}$. Since $w_0$ was an arbitrary point in $H\mcg^{-1}(0)\cap \tilde{Z}$, we get \eqref{eq_hyp7}. Together with \eqref{eq_hyp6} this gives us \eqref{eq_hyp5}, which we can combine with \eqref{eq_hyp4} to finally get \eqref{eq_cl}, i.e. $\overline{S\mcg}=H\mcg^{-1}(0)$.

Now that we have found a nice formula for $\overline{S\mcg}$, we can prove that it is a hypersurface in $\widetilde{X}$. For this we first compute the partial derivative of $H\mcg$ with respect to $P_r$:
\[\frac{\partial H\mcg}{\partial P_r}=2\psi^TB\psi.\]
We can again use that $2\psi^TB\psi$ is always positive. Hence $H\mcg$ has no singular points and in particular 0 is a regular value. Thus we get that $\overline{S\mcg}=H\mcg^{-1}(0)$ is a (smooth) hypersurface in $\widetilde{X}$. Since $\Pi|_{\widetilde{X}}$ is a diffeomorphism onto its image, therefore $\Pi(H\mcg^{-1}(0))=S_b\cap \Pi(\widetilde{X})$ is a smooth hypersurface in $\Pi(\widetilde{X})$, which finishes the proof.
\end{proof}

Now we are finally ready to prove the main theorem of this paper:

\vspace*{0.25cm}
\noindent\textit{Proof of Theorem \ref{thm_main}:} By Lemma \ref{lem_Sb} we know that the $b$-image of $S$ in $(X,Z)$ is a smooth hypersurface $S_b$, which is described as zero level set of two different Hamiltonians:\linebreak $S_b\cap \Pi(M\times T^*\R^m\times \{-1,1\}=\Pi(H^{-1}(0)\times\{-1,1\})$, whereas $S_b \cap \Pi(\widetilde{X})= \Pi\left(H\mcg^{-1}(0)\right)$. By Lemma  \ref{lem_tlvf}, to show that $S_b$ is a contact type hypersurface in $(X,Z, \omega \oplus \omega_b)$ it suffices to find a Liouville vector field that is transverse to $S_b$.

Let $Y$ be the global Liouville vector field on $(M, \omega)$ as in Axiom \eqref{defH0} and let $V$ be the Liouville vector field defined in \eqref{eq_defV}. We will show that $Y +V$ is transverse to $S_b$.

Note that $\{M \times \pi (T^*\R^k \times \{-1,1\}), M\times \pi( \{|z|<\sqrt{2}\})\}$ is an open cover of $X$. Therefore, it suffices to show that 
\begin{align}
dH\left(Y+\widehat{V}\right)\neq 0 & \quad \text{on} \quad S \subseteq M\times T^*\R^k, \label{eq_dHVhat}\\
dH\mcg\left(Y+\widetilde{V}\right)\neq 0 & \quad   \text{on} \quad H\mcg^{-1}(0) \cap \left(M \times \{|z|<\sqrt{2}\}\right)\subseteq \widetilde{X},\label{eq_dHVtilde}
\end{align}
where the last equality comes from the fact that $\Pi(H\mcg^{-1}(0))=S_b\cap \Pi(\widetilde{X})$.

Using \eqref{eq_defV} we can calculate:
\begin{align*}
dH\left(Y+\widehat{V}\right)& = \left(dH_0+p^TBdq+qB^Tdp\right)\left( Y+(p+h(|q|)q)\partial_p\right)\\
& = dH_0(Y)(x)+p^TBq+h(|q|)q^TBq\\
& = H(x, q, p)+(dH_0(Y)-H_0)(x)+\frac{1}{2}h(|q|)q^T(B+B^T)q.
\end{align*}
By definition $H$ vanishes on $S$ so the first summand vanishes. On the other hand, we know that by \eqref{defH0} the second summand is strictly positive. Finally, by assumption the matrix $B+B^T$ is positive definite and $h$ is non-negative, so the third summand is non-negative. Consequently, the sum is positive on $S:= H^{-1}(0)$, which proves \eqref{eq_dHVhat}.

To prove \eqref{eq_dHVtilde} first observe that
\[
\widetilde{V}(z,P_r,\phi,\eta)= (P_r+1)\partial_{P_r}+\eta\partial_\eta \quad \text{on} \quad \left\lbrace |z|<\sqrt{2}\right\rbrace \subseteq \widetilde{W}.
\]
Therefore, using \eqref{eq_Hmcg} we obtain
\begin{align*}
dH\mcg\left(Y+\widetilde{V}\right) & = \left( z^2 dH_0+ 2\psi^T B \psi dP_r+z^2\psi B\Psi d\eta\right)\left( Y+(P_r+1)\partial_{P_r}+\eta\partial_\eta\right)\\
& = z^2 dH_0(Y)+ 2 (P_r +1)\psi^T B \psi+z^2\psi B\Psi \eta\\
& = H\mcg +z^2 (dH_0(Y)-H_0) + \psi^T (B+B^T) \psi.
\end{align*}
Recall, that by assumption \eqref{defH0} $dH_0(Y)-H_0$ is strictly positive. On the other hand, by assumption $B+B^T$ is positive definite and $|\psi|=1$, thus 
\[dH\mcg\left(Y+\widetilde{V}\right)\geq \psi^T (B+B^T) \psi>0 \quad \text{on} \quad H\mcg^{-1}(0),\]
which proves \eqref{eq_dHVtilde}. Combining \eqref{eq_dHVhat} and \eqref{eq_dHVtilde} we conclude that the Liouville vector field $Y+ V$ is transverse to $S_b$, which means that the $1$-form $\iota_{Y}\omega+ \iota_V\omega_b$ restricted to $TS_b$ defines a $b$-contact structure.

\hfill $\square$

\appendix
\section{Non-uniqueness of \texorpdfstring{$b$}{b}-symplectic structures}\label{app_nonun}
In these appendices we will present some counterexamples to give an intuition of what can go wrong in the setting of $b$-compatible manifolds. First, we will investigate the uniqueness of $b^m$-symplectic structures. In particular, it is a natural question to ask whether they are fully determined by their behaviour away from the critical set. This can be formulated more precisely as:
\begin{question}
Let $\omega_1$ and $\omega_2$ be two $b^m$-symplectic forms on a $b$-manifold $(M,Z)$ such that the symplectic manifolds $(M\setminus Z,\omega_1)$ and $(M\setminus Z,\omega_2)$ are symplectomorphic. Do the two\linebreak $b^m$-symplectic manifolds $(M,Z,\omega_1)$ and $(M,Z,\omega_2)$ necessarily have to be $b$-symplectomorphic?
\end{question}
As it turns out the answer to this question is no. Consider the manifold
\begin{equation*}
M:=T^*\T^2=\left\{((\phi,\theta,s,t)\,\big|\,\phi,\theta\in\R/\pi\Z;s,t\in\R\right\}
\end{equation*}
and the hypersurface $Z:=\{\phi=0\}\subset M$. Then we can define on the $b$-manifold $(M,Z)$ the following two $b^2$-symplectic forms:
$$
\omega_1:=\frac{1}{\sin{\phi}^2}d\phi\wedge d\theta+ds\wedge dt,
\qquad\omega_2:=\frac{1}{\sin{\phi}^2}d\phi\wedge ds+d\theta\wedge dt.
$$
It is easy to see that $\omega_1$ and $\omega_2$ are well-defined, closed and non-degenerate $b^2$-2-forms on $(M,Z)$. Hence $(M,Z,\omega_1)$ and $(M,Z,\omega_2)$ are indeed well-defined $b^2$-symplectic manifolds. Furthermore, one can show that they are not $b$-symplectomorphic. For this we simply notice that $\omega_2$ is exact while $\omega_1$ is not (see Lemma \ref{lem_exact}). This then directly implies that they are not $b$-symplectomorphic since any $b$-symplectomorphism induces an isomorphism on the \linebreak$b^2$-cohomology. For our purposes it is enough to know that $b^2$-cohomology is well-defined, for details see \cite{scott2016} for the original construction. Finally, in Lemma \ref{lem:symp} we show that after removing the singular set $Z$ we obtain two symplectomorphic manifolds $(M\setminus Z,\omega_1)$ and $(M\setminus Z,\omega_2)$.
\begin{remark}
Note that $\omega_1$ and $\omega_2$ were defined using the same coordinates, so they automatically live in the same space of $b^2$-forms and we do not have to worry about specifying the jet-data mentioned in Remark \ref{rem_jet}.
\end{remark}
\begin{lemma}\label{lem_exact}
The $b^2$-2-form $\omega_2$ is exact while $\omega_1$ is not.
\end{lemma}
\begin{proof}
Indeed, we have
\begin{equation*}
\omega_2=d\left(-\frac{s}{\sin{\phi}^2}d\phi-td\theta\right).
\end{equation*}

Now let us assume that $\omega_1$ is also exact, i.e. there exists a $b^2$-1-form $\alpha$ with $\omega_1=d\alpha$. Then we get two functions defined by $\gamma:=\alpha(\partial_\phi)$ and $\xi:=\alpha(\partial_\theta)$ respectively. Since $\partial_\theta$ is a $b^2$-vector field while $\partial_\phi$ is not, $\xi$ is a well-defined smooth function on all of $M$ whereas $\gamma$ only when restricted to $M\setminus Z$. For these two functions we get the equation:
\begin{equation*}
\frac{1}{\sin{\phi}^2}=\omega_1(\partial_\phi,\partial_\theta)=d\alpha(\partial_\phi,\partial_\theta)=\mathcal{L}_{\partial_\phi}\left(\alpha(\partial_\theta)\right)-\mathcal{L}_{\partial_\theta}\left(\alpha(\partial_\phi)\right)-\alpha\left([\partial_\phi,\partial_\theta]\right)=\frac{\partial\xi}{\partial \phi}-\frac{\partial\gamma}{\partial \theta},
\end{equation*}
where we used the fact that $[\partial_\phi,\partial_\theta]=0$. If we fix some $\phi\neq0$ and $s=t=0$ we can integrate over $\theta\in[0,\pi]$ to get:
\begin{equation*}
\frac{\pi}{\sin{\phi}^2}=\int_0^{\pi}\frac{\partial \xi}{\partial \phi}d\theta-\int_0^{\pi}\frac{\partial \gamma}{\partial \theta}d\theta=\int_0^{\pi}\frac{\partial\xi}{\partial \phi}d\theta,
\end{equation*}
where we used that $\gamma$ has to be $\pi$-periodic in $\theta$. But since $\xi$ is a smooth function on all of $M$, we can bound all of its partial derivatives on the compact set $\{s=t=0\}$. Thus this equality cannot be satisfied for $\phi$ approaching $0$ and we reach a contradiction.
\end{proof}
However, we can now show that if we remove the singular hypersurface then the two resulting symplectic manifolds are symplectomorphic. Hence we get that in general a \linebreak$b^m$-symplectic structure is not uniquely determined by the symplectic structure away from the singular set.
\begin{lemma}\label{lem:symp}
The two symplectic manifolds $(M\setminus Z,\omega_1)$ and $(M\setminus Z,\omega_2)$ are symplectomorphic.
\end{lemma}
\begin{proof}
Recall, that the cotangent function $\cot: (0,\pi) \to \R$ is a diffeomorphism. Thus we can define a new diffeomorphism on all of $M\setminus Z$:
\begin{align*}
f:M\setminus Z&=(0,\pi)\times\R/\pi\Z\times\R^2\to M\setminus Z,\\
f(\phi,\theta,s,t)&:=\left(\cot^{-1}(t),\theta,s,\cot(\phi)\right).
\end{align*}
Then we calculate $f^*(dt)=\frac{-1}{\sin{\phi}^2}d\phi$ and similarly $f^*(\frac{1}{\sin{\phi}^2}d\phi)=-dt$ with which we just get:
\begin{equation*}
f^*(\omega_1)=f^*\left(\frac{1}{\sin{\phi}^2}d\phi\wedge d\theta+ds\wedge dt\right)=-dt\wedge d\theta+ds\wedge\frac{-1}{\sin{\phi}^2}d\phi=\omega_2.
\end{equation*}
Hence $f$ is a symplectomorphism between $(M\setminus Z,\omega_1)$ and $(M\setminus Z,\omega_2)$ as required.
\end{proof}
\section{Properties of \texorpdfstring{$b$}{b}-images of hypersurfaces}
In this appendix we will explore the concept of the $b$-image, as defined in Definition \ref{def_img}. First of all, one might wonder whether the $b$-image retains submanifold structures, or in particular:
\begin{question}
Let $S$ be a hypersurfaces in $M$. Is its $b$-image $S_b$ in any $b$-compatible manifold always also a hypersurface?
\end{question}
Unfortunately, this is not true in general. In fact, for $M=T^\ast\R^n$ we can write down an explicit counterexample by taking the $b$-manifold $(W,Z)$ as defined in Section \ref{sec_glm} as our $b$-compatible manifold. To avoid confusion we will keep the same notation.
\begin{lemma}\label{lem_hypcou2}
The hypersurface $U:=\{q_1+\sum_{i=2}^nq_i^2=\frac{1}{4}\}$ is a hypersurface in $T^\ast\R^n$, but $U_b$ is not a hypersurface in $W$.
\end{lemma}
\begin{proof}
By construction, $U$ is given by $H^{-1}(\frac{1}{4})$ for the smooth function:
\[H=q_1+\sum_{i=2}^nq_i^2,\]
for which $\frac{1}{4}$ is a regular value. Next recall the spherical change of coordinates\linebreak $F:\R^n\setminus\{0\}\to \R_+\times \s^{n-1}=:Q$ as defined in \eqref{F}. Using this map we can write the equation for $U$ in $n$-spherical coordinates by defining $\psi_i:=q_i/r=\frac{q_i}{\sqrt{\sum_{i=1}^nq_i^2}}$:
\begin{equation*}
\sum_{i=2}^nq_i^2=r^2-q_1^2=r^2(1-\psi_1^2).
\end{equation*}
With this we get that $U$ as a subset of $T^\ast Q$ is given by:
\begin{equation}\label{eq_shy}
U=\left\{r\psi_1+r^2(1-\psi_1^2)=\frac{1}{4}\right\}.
\end{equation}
Now we can start to look at $U_b$. Recall the definition of the McGehee transform $\tau\mcg$ from \eqref{eq_tau}. By substituting $r=2/z^2$ into \eqref{eq_shy} we get:
\[\tau\mcg^{-1}(U)=\left\{2z^2\psi_1+4(1-\psi_1^2)=\frac{1}{4}z^4,z\neq 0\right\}.\]
Solving for $z$ gives us $\{z=\pm 2\sqrt{\psi_1\pm 1},z\neq 0\}$. Since $\psi_1$ is a spherical coordinate we have $\psi_1\in [-1,1]$. Consequently, the equation $z=\pm 2\sqrt{\psi_1-1}$ has no non-zero solutions for $z$, so in fact we have:
\[\tau\mcg^{-1}(U)=\{z=\pm 2\sqrt{\psi_1+1},z\neq 0\}.\]
Since $z\to 0$ is equivalent to $\psi_1\to-1$ on this set, the closure of this is simply:
\[U_b=\{z=\pm 2\sqrt{\psi_1+1}\}\m{.}\]
This set has a singularity at $\{\psi_1=-1,z=0\}$, which we will verify by computing its tangent space, thus proving that $U_b$ is not a hypersurface.

Let us now show that $U_b$ has a singularity at $\{\psi_1=-1,z=0\}$. We will prove it by contradiction assuming that $U_b$ is a hypersurface in $W$. Since $U$ does not intersect the set $\{r=0\}$ its $b$-image $U_b$ is entirely contained in an open neighbourhood of $Z$, which we can identify with $\R^2\times T^*\s^{n-1}$. The space $T^*\s^{n-1}$ itself can be seen as a submanifold of codimension 2 in $\R^{2n}$ defined by the equations $\sum_{i=1}^n\psi_i^2=1$ and $\sum_{i=1}^n\psi_i\zeta_i=0$ (where we use coordinates $(\psi_1,\dots,\psi_n,\zeta_1,\dots,\zeta_n)$ for $\R^{2n}$). Thus we can see $U_b$ as a codimension 3 submanifold of the euclidean space $\R^{2n+2}$  (where $z$ and $P_r$ are the two additional coordinates). Furthermore, we know that it is given by the three equations:
\begin{equation*}
\sum_{i=1}^n\psi_i^2=1,\quad\sum_{i=1}^n\psi_i\zeta_i=0\quad\text{and}\quad \frac{z^2}{4}=\psi_1+1.
\end{equation*}
This gives us that the tangent space of $U_b$ satisfies:
\begin{equation*}
TU_b\subseteq \operatorname{ker}\left(\sum_{i=1}^n\psi_id\psi_i\right)\cap \operatorname{ker}\left(\sum_{i=1}^n\zeta_id\psi_i+\sum_{i=1}^n\psi_id\zeta_i\right)\cap \operatorname{ker}\left(\frac{z}{2}dz-d\psi_1\right).
\end{equation*}
Note that these 1-forms are linearly dependent on $U_b$ if and only if $z=0$. Hence this is an equality whenever $z\neq 0$. To compute $T^*U_b\big|_{\{z=0\}}$ we will need to use a continuity argument.

If we restrict ourselves to the set
\begin{equation*}
A^+:=\{\psi_1<0 <\psi_2, z,\ \psi_3=\dots=\psi_n=\zeta_3=\dots=\zeta_n=0\}
\end{equation*}
we can solve two of the above equations for $\psi_1$:
\begin{equation*}
\psi_2=\sqrt{1-\psi_1^2}\quad\text{and}\quad z=2\sqrt{\psi_1+1}.
\end{equation*}
We can substitute this into the above formula for $TU_b$ to get:
\begin{align*}
TU_b\big|_{A^+}=\operatorname{ker}\left(\psi_1d\psi_1+\sqrt{1-\psi_1^2}d\psi_2\right)&\cap \operatorname{ker}\left(\zeta_1d\psi_1+\zeta_2d\psi_2+\psi_1d\zeta_1+\sqrt{1-\psi_1^2}d\zeta_2\right)
\\&\cap \operatorname{ker}\left(\sqrt{\psi_1+1}dz-d\psi_1\right).
\end{align*}
Thus we have that the tangent space on $A^+$ contains the vector fields $\partial_{P_r}$, $\partial_{\psi_i}, \partial_{\zeta_i}$ for $i>2$ and the two vector fields
\begin{align*}
X_1^+ & :=\sqrt{1-\psi_1^2}\partial_{\zeta_1}-\psi_1\partial_{\zeta_2},\\
X_2^+ &:=\partial_z + \sqrt{\psi_1+1}\partial_{\psi_1}-\frac{\psi_1}{\sqrt{1-\psi_1}}\partial_{\psi_2}+\left(\frac{-\zeta_1\sqrt{\psi_1+1}}{\psi_1}+\frac{\zeta_2}{\sqrt{1-\psi_1}}\right)\partial_{\zeta_1}.
\end{align*}
In the limit we have
$$
\lim_{\psi_1\to -1}X_1^+ = \partial_{\zeta_2}\qquad \textrm{and}\qquad \lim_{\psi_1\to -1}X_2^+=\partial_z+\frac{1}{\sqrt{2}}\partial_{\psi_2}+\frac{\zeta_2}{\sqrt{2}}\partial_{\zeta_1}.
$$
Next we can repeat this calculation for the set
\begin{equation*}
A^-:=\{\psi_1,\psi_2<0<z,\ \psi_3=\dots=\psi_n=\zeta_3=\dots=\zeta_n=0\}.
\end{equation*}
Similarly, this gives us that $TU_b\big|_{A^-}$ includes vector fields $\partial_{P_r}$, $\partial_{\psi_i}, \partial_{\zeta_i}$ for $i>2$ and the two vector fields
\begin{align*}
X_1^- & :=\sqrt{1-\psi_1^2}\partial_{\zeta_1}+\psi_1\partial_{\zeta_2},\\
X_2^- & :=\partial_z + \sqrt{\psi_1+1}\partial_{\psi_1} +\frac{\psi_1}{\sqrt{1-\psi_1}}\partial_{\psi_2} - \left(\frac{\zeta_1\sqrt{\psi_1+1}}{\psi_1}+\frac{\zeta_2}{\sqrt{1-\psi_1}}\right)\partial_{\zeta_1},
\end{align*}
which only differ by signs from the ones calculated above. In the limit we have
$$
\lim_{\psi_1\to -1}X_1^- = \partial_{\zeta_2}\qquad \textrm{and}\qquad \lim_{\psi_1\to -1}X_2^-=\partial_z-\frac{1}{\sqrt{2}}\partial_{\psi_2}-\frac{\zeta_2}{\sqrt{2}}\partial_{\zeta_1}.
$$
By continuity, the tangent space of $U_b$ at $\{\psi_1=-1,\psi_2=\dots=\psi_n=z=0\}$ then needs to contain all the $2n-2$ vector fields $\{\partial_{\psi_3},\dots,\partial_{\psi_n},\partial_{P_r},\partial_{\zeta_2},\dots,\partial_{\zeta_n}\}$ and the limits of both $X_2^+$ and $X_2^-$. However, these all are clearly linearly independent so they cannot all lie in the same $2n-1$ dimensional subspace. Hence by contradiction $U_b$ cannot be a hypersurface in $W$.
\end{proof}
Next, it is a natural question to ask whether the definition of the $b$-image depends on the diffeomorphisms used to define $b$-compatibility. Or equivalently, one could ask the following question:
\begin{question}
Given two diffeomorphic hypersurfaces $S$ and $U$ in a manifold $M$, are their $b$-images $S_b$ and $U_b$ also diffeomorphic in any $b$-compatible manifold?
\end{question}
Again, the answer is no. Furthermore, we can also find counterexamples defined in $(W,Z)$.
\begin{lemma}
The hypersurfaces $S:=\{q_1=\frac{1}{4}\}$ and $U:=\{q_1+\sum_{i=2}^nq_i^2=\frac{1}{4}\}$ in $T^\ast\R^n$ are diffeomorphic, but their $b$-images $S_b$ and $U_b$ are not diffeomorphic in $(W,Z)$.
\end{lemma}
\begin{proof}
The two hypersurfaces $\{q_1=\frac{1}{4}\}$ and $\{q_1+\sum_{i=2}^nq_i^2=\frac{1}{4}\}$ in $\R^n$ are diffeomorphic via the map:
\[(q_1,q_2,\dots,q_n)\mapsto\left(q_1-\sum_{i=2}^nq_i^2,q_2,\dots,q_n\right).\]
Next we can lift this to get a diffeomorphism between $S$ and $U$ by looking at the induced map on the cotangent bundle. From Lemma \ref{lem_hypcou2} we know that $U_b$ isn't a hypersurface in $(W,Z)$. However, we will show that $S_b$ is a hypersurface in $(W,Z)$.

Similarly as with $U$, the hypersurface $S$ also does not intersect the set $\{r=0\}$, hence we can also see it as a hypersurface in $T^\ast Q$. By substituting $q_1=r\psi_1$ as above we get that $S$ in $T^\ast Q$ is given by $S=\{r \psi_1=1\}$. When applying the McGehee transform and substituting $r=2/z^2$ into this we get $\{2\psi_1=z^2\}$. We will show that this is in fact exactly $S_b$ and it does define a hypersurface in $W$.

More precisely, similar as in the proof of the previous lemma, we can see $S_b$ as the subset of $\R^{2+2n}$ defined by the three equations:
\begin{equation*}
f:=\sum_{i=1}^n\psi_i^2-1=0\text{,}\quad g:=\sum_{i=1}^n\psi_i\zeta_i=0\quad\text{and}\quad h:=2\psi_1-z^2=0.
\end{equation*}

Then we notice that the 1-forms $df$, $dg$ and $dh$ are linearly dependent only on the sets 
$$
B:=\{\psi_1=\dots=\psi_n=0\}\qquad \textrm{and} \qquad C:=\{\psi_2=\dots=\psi_n=z=0\}.
$$
However, $f\big|_B$ is always equal to $-1$. Hence the intersection $f^{-1}(0)\cap B$ is empty. Similarly, $h\big|_C$ is just equal to $2\psi_1$ and so the intersection $h^{-1}(0)\cap C$ is contained in $B$. Thus the intersection $f^{-1}(0)\cap h^{-1}(0)\cap C$ is empty as above. Combined, this gives us that $df$, $dg$ and $dh$ are linearly independent on $S_b=f^{-1}(0)\cap g^{-1}(0)\cap h^{-1}(0)$. This in turn means that $S_b$ is a smooth hypersurface in $f^{-1}(0)$ and hence also in $W$.

In the end, since we have that $S_b$ is a hypersurface in $W$ and $U_b$ is not, they clearly cannot be diffeomorphic.
\end{proof}
\begin{remark}
Note that the lifted diffeomorphism between $S$ and $U$ is in fact a symplectomorphism by the properties of the canonical symplectic structure on the cotangent bundle. Hence we immediately also get that the properties of a $b$-image in a $b$-symplectically compatible manifold heavily depend on the chosen symplectomorphisms.
%In particular this also means that they aren't $b^3$-symplectomorphic in $W$. In fact, they are not even homeomorphic. This means that 
\end{remark}

\printbibliography

\end{document}